\documentclass{gtart_h}  

\def\ifplaintex{\expandafter\ifx\csname documentclass\endcsname\relax}

\ifplaintex 
\else
\fi

\expandafter\ifx\csname beginpicture\endcsname\relax
\expandafter\ifx\csname documentclass\endcsname\relax
\input pictex \else
\input prepictex \input pictex \input postpictex \fi\fi

\def\gt{{\mathsurround=0pt\it $\cal G\mskip-2mu$eometry \&\ 
$\cal T\!\!$opology}}        

\def\gtp{{\mathsurround=0pt\it $\cal G\mskip-2mu$eometry \&\ 
$\cal T\!\!$opology $\cal P\!$ublications}}  

\def\lognumber#1{\def\thelognumber{#1}}
\def\volumenumber#1{\def\thevolumenumber{#1}}
\def\papernumber#1{\def\thepapernumber{#1}}
\def\volumeyear#1{\def\thevolumeyear{#1}}

\def\pagenumbers#1#2{\def\startpage{#1}\def\finishpage{#2}}
\def\published#1{\def\publishdate{#1}}
\def\proposed#1{\def\theproposer{#1}}
\def\seconded#1{\def\theseconders{#1}}
\def\received#1{\def\receiveddate{#1}}

\def\accepted#1{\def\accepteddate{#1}}
\def\asciititle#1{\def\theasciititle{#1}}

\long\def\asciiabstract#1{\long\def\theasciiabstract{#1}}
\def\asciikeywords#1{\def\theasciikeywords{#1}}

\def\shorttitle#1{\def\theshorttitle{#1}}

\let\\\par\let\thelognumber\relax
\let\thevolumenumber\relax\let\thepapernumber\relax
\let\thevolumeyear\relax\let\thesamplenumber\relax\let\startpage\relax
\let\finishpage\relax\let\publishdate\relax\let\receiveddate\relax
\let\reviseddate\relax\let\accepteddate\relax\let\theasciititle\relax
\let\theasciiauthors\relax
\let\theasciiabstract\relax\let\theasciikeywords\relax
\let\theasciiemail\relax\let\theshortauthors\relax\let\theshorttitle\relax

\long\def\maketitlep{   

\count0=\startpage

\gt\hfill      
\beginpicture
\setcoordinatesystem units <0.33truein, 0.33truein> point at 2.2 0.9
\setplotsymbol ({$\cal G$})
\plotsymbolspacing=9truept
\circulararc 315 degrees from 0 1 center at 0 0
\setplotsymbol ({$\cal T$})
\circulararc 315 degrees from 1 -1 center at 1 0
\endpicture
\break
{\small\ifx\thesamplenumber\relax 
Volume \else Sample
\fi\thevolumenumber\ (\thevolumeyear)
\startpage--\finishpage\nl
Published: \publishdate}
\vglue 0.5truein plus 0.4fil minus 0.1truein

{\parskip=0pt\leftskip 0pt plus 1fil\def\\{\par\smallskip}{\ifplaintex\large
\else\Large\fi\bf\thetitle}\par\medskip}   

\vglue 0pt plus 0.1fil 

{\parskip=0pt\leftskip 0pt plus 1fil\def\\{\par}{\sc\theauthors}
\par\medskip}

\vglue 0pt plus 0.1fil 

{\small\parskip=0pt\let\newline\\
{\leftskip 0pt plus 1fil\def\\{\par}{\sl\theaddress}\par}
\expandafter\ifx\theemail\relax    
\relax\else\vglue 5pt plus 0.02fil minus 2pt\def\\{\stdspace{\rm 
and}\stdspace} 
\cl{Email:\stdspace\tt\theemail}\fi
\ifx\theurl\relax                  
\relax\else\vglue 5pt plus 0.02fil minus 2pt\def\\{\stdspace{\rm 
and}\stdspace}
\cl{URL:\stdspace\tt\theurl}\fi\par}

\vglue 7pt plus 0.3fil minus 3pt

{\bf Abstract}
\vglue 5pt plus 0.1fil minus 2pt

\theabstract

\vglue 7pt plus 0.3fil minus 3pt

{\bf AMS Classification numbers}\quad Primary:\quad \theprimaryclass

Secondary:\quad \thesecondaryclass

\vglue 5pt plus 0.3fil minus 2pt

{\bf Keywords:}\quad \thekeywords

\vglue 10pt plus 0.5fil minus 5pt

{\small  Proposed: \theproposer\hfill Received: \receiveddate\nl
Seconded: \theseconders\hfill 
\ifx\reviseddate\relax                         
Accepted: \accepteddate                        
\else
Revised: \reviseddate                          
\fi}
\eject
}       

\let\maketitlepage\maketitlep
\let\maketitle\maketitlepage

\font\phead=cmsl9 scaled 950
\font=cmsl9 scaled 1050
\font\pnum=cmbx10 scaled 913
\font\lnum=cmbx10 
\font\pfoot=cmsl9 scaled 950
\font=cmsl9 scaled 1050
\ifplaintex
\headline{\vbox to 0pt{\vskip -4.5mm\line{\small\phead\ifnum
\count0=\startpage ISSN 1364-0380 (on line)
1465-3060 (printed) \hfill {\pnum\folio}\else\ifodd\count0\def\\{ }%
\ifx\theshorttitle\relax\thetitle\else\theshorttitle\fi\hfill{\pnum\folio}
\else\def\\{ and }{\pnum\folio}\hfill\ifx\theshortauthors\relax\theauthors
\else\theshortauthors\fi\fi\fi}\vss}}
\footline{\vbox to 0pt{\vglue 0mm\line{\small\pfoot\ifnum\count0=\startpage
\copyright\ \gtp\hfill\else
\gt, Volume \thevolumenumber\ (\thevolumeyear)\hfill\fi}\vss
}}
\else
\makeatletter
\def\@oddhead{{\small\lhead\ifnum\count0=\startpage ISSN 1364-0380 (on line)
1465-3060 (printed) \hfill {\lnum\number\count0}\else\ifodd\count0
\def\\{ }\ifx\theshorttitle\relax \thetitle \else\theshorttitle\fi\hfill
{\lnum\number\count0}\else\def\\{ and }{\lnum\number\count0}
\hfill\ifx\theshortauthors\relax 
\theauthors\else\theshortauthors\fi\fi\fi}}\def\@evenhead{\@oddhead}
\def\@oddfoot{\small\lfoot\ifnum\count0=\startpage\copyright\ \gtp\hfill\else
\gt, Volume \thevolumenumber\ (\thevolumeyear)\hfill\fi}
\def\@evenfoot{\@oddfoot}
\makeatother
\fi

\newwrite\gtoutfile
\long\gdef\makeheadfile{  
{\def\\{, }\def\s{ }
\immediate\openout\gtoutfile head.xxx
\immediate\write\gtoutfile{Proxy-for: \ifx\theasciiauthors\relax
\theauthors\else\theasciiauthors\fi\s<\ifx\theasciiemail\relax\theemail\else\theasciiemail\fi>}
\immediate\write\gtoutfile{\noexpand\\}
\immediate\write\gtoutfile{Authors: \ifx\theasciiauthors\relax
\theauthors\else\theasciiauthors\fi}
{\def\\{ }\immediate\write\gtoutfile{Title: \ifx\theasciititle\relax
\thetitle\else\theasciititle\fi}}
\immediate\write\gtoutfile{Subj-class: GT or SG or MG etc}
\immediate\write\gtoutfile{MSC-class: \theprimaryclass\ifx\thesecondaryclass\relax\else, \thesecondaryclass\fi}
\immediate\write\gtoutfile{Journal-ref: Geom. Topol. \thevolumenumber
(\thevolumeyear) \startpage-\finishpage}
\immediate\write\gtoutfile{Comments: Published by Geometry and Topology at}
\immediate\write\gtoutfile{\s\s http://www.maths.warwick.ac.uk/gt/GTVol\thevolumenumber/paper\thepapernumber.abs.html}
\immediate\write\gtoutfile{\noexpand\\}
\immediate\write\gtoutfile{}
\ifx\theasciiabstract\relax
\immediate\write\gtoutfile{\theabstract}\else
\immediate\write\gtoutfile{\theasciiabstract}\fi
\immediate\write\gtoutfile{}
\immediate\write\gtoutfile{\noexpand\\}
\immediate\write\gtoutfile{}
\immediate\closeout\gtoutfile}}  

\def\maketitlepage{\maketitlep\makeheadfile}
\let\maketitle\maketitlepage

\lognumber{398}
\volumenumber{8}\papernumber{14}\volumeyear{2004}
\pagenumbers{565}{610}
\received{18 December 2003}
\published{31 March 2004}
\accepted{26 March 2004}
\proposed{Yasha Eliashberg}
\seconded{Ronald Fintushel, Ronald Stern}

\usepackage{amssymb,amsmath}
\newtheorem{prop}{Proposition}
\newtheorem{lemma}[prop]{Lemma}
\newtheorem{thm}[prop]{Theorem}
\newtheorem{cor}[prop]{Corollary}

\DeclareMathOperator{\cok}{coker} \DeclareMathOperator{\ind}{ind}
\DeclareMathOperator{\ep}{\epsilon}
\DeclareMathOperator{\dbar}{\bar{\partial}}
\numberwithin{equation}{section} \numberwithin{prop}{section}

\title{The Gromov invariant and the Donaldson--Smith\\standard surface count}
\asciititle{The Gromov invariant and the Donaldson-Smith standard surface count}
\shorttitle{Gromov invariant and Donaldson--Smith surface count}
\author{Michael Usher}
\address{Department of Mathematics, MIT\\Cambridge, MA
02139--4307, USA} 
\email{usher@math.mit.edu}
\begin{document}

\begin{abstract} Simon Donaldson and Ivan Smith recently studied
symplectic surfaces in symplectic 4--manifolds $X$ by introducing
an invariant $\mathcal{DS}$ associated to any Lefschetz fibration
on blowups of $X$ which counts holomorphic sections of a relative
Hilbert scheme that is constructed from the fibration.  Smith has
shown that $\mathcal{DS}$ satisfies a duality relation identical
to that satisfied by the Gromov invariant $Gr$ introduced by
Clifford Taubes, which led Smith to conjecture that
$\mathcal{DS}=Gr$ provided that the fibration has high enough
degree. This paper proves that conjecture. The crucial technical
ingredient is an argument which allows us to work with curves $C$
in the blown-up 4--manifold that are made holomorphic by an almost
complex structure which is integrable near $C$ and with respect to
which the fibration is a pseudoholomorphic map.
\end{abstract}

\asciiabstract{%
Simon Donaldson and Ivan Smith recently studied symplectic surfaces in
symplectic 4-manifolds X by introducing an invariant DS associated to
any Lefschetz fibration on blowups of X which counts holomorphic
sections of a relative Hilbert scheme that is constructed from the
fibration.  Smith has shown that DS satisfies a duality relation
identical to that satisfied by the Gromov invariant Gr introduced by
Clifford Taubes, which led Smith to conjecture that DS=Gr provided
that the fibration has high enough degree. This paper proves that
conjecture. The crucial technical ingredient is an argument which
allows us to work with curves C in the blown-up 4-manifold that are
made holomorphic by an almost complex structure which is integrable
near C and with respect to which the fibration is a pseudoholomorphic
map.}


\asciikeywords{Pseudoholomorphic curves, symplectic Lefschetz
fibrations, Gromov-Witten invariants}

\maketitlepage

\section{Introduction}
Let $(X,\omega)$ be a symplectic 4--manifold.  Since the
publication of Simon Donaldson's famous paper \cite{Don} it has
been realized that a fruitful way of studying $X$ is to construct
a symplectic Lefschetz fibration $f\co X'\to S^2$ on a suitable
blow-up $X'$ of $X$.  One application of Lefschetz fibration
techniques has been the work of Donaldson and Ivan Smith in
\cite{DS} and \cite{Smith} toward re-proving results concerning
holomorphic curves in $X$ which were originally obtained by Cliff
Taubes in his seminal study of the Seiberg--Witten equations on
symplectic manifolds.  In \cite{Taubes}, Taubes constructs a
``Gromov invariant'' $Gr(\alpha)$ which counts embedded, not
necessarily connected, pseudoholomorphic submanifolds of $X$ which
are Poincar\'{e} dual to a class $\alpha\in{}H^{2}(X;\mathbb{Z})$,
and in his other papers (collected in \cite{book}) he identifies
$Gr$ with the Seiberg--Witten invariants.  From the
charge--conjugation symmetry in Seiberg--Witten theory there then
follows the surprising Taubes duality relation that, where
$\kappa$ is the canonical class of $X$ (ie, the first
Chern class of the cotangent bundle), $Gr(\alpha)=\pm Gr(\kappa-
\alpha)$, provided that $b^{+}(X)>1$.

One might reasonably expect that a formula such as the Taubes
duality relation could be proven in a more hands-on way than that
provided by Seiberg--Witten theory, and Donaldson and Smith have
indeed provided a somewhat more intuitive framework for
understanding it.  After perturbing $\omega$ to make its
cohomology class rational and then scaling it to make it integral,
Donaldson's construction gives, for large enough $k$, symplectic
Lefschetz pencils $f_k\co  X\setminus B_k \to S^2$ ($B_k$ being a
set of $k^2 [\omega]^2$ points obtained as the common vanishing
locus of two sections of a line bundle over $X$) which lift to
symplectic Lefschetz fibrations $f'_k\co  X'_k \to S^2$ where
$\pi_k\co  X'_k \to X$ is the blowup of $X$ along $B_k$; the
fibers of $f'_k$ are Poincar\'e dual to $k\pi_k ^{*}[\omega]$.
 From any symplectic Lefschetz fibration $f\co X'\to S^2$ and for any natural number $r$ Donaldson
and Smith \cite{DS} construct the ``relative Hilbert scheme''
$F\co X_{r}(f)\to S^2$ whose fiber over a regular value $t$ of $f$
is the symmetric product $S^{r}f^{-1}(t)$; this is a smooth
manifold that can be given a (continuous family of) symplectic
structure(s) by the Thurston trick.  A section of $F$ then
naturally corresponds to a closed set in $X'$ which intersects
each fiber of $f$ $r$ times (possibly counting multiplicities). So
if we take an almost complex structure $j$ on $X'$ with respect to
which the fibration $f\co X'\to S^2$ is a pseudoholomorphic map
(so that in particular the fibers of $f$ are $j$--holomorphic and
therefore intersect other $j$--holomorphic curves locally
positively), then a holomorphic curve Poincar\'e dual to some
class $\alpha$ and not having any fiber components will, to use
Smith's words, ``tautologically correspond'' to a section of $X_r
(f)$. This section will further be holomorphic with respect to the
almost complex structure $\mathbb{J}_j$ on $X_{r}(f)$ obtained
from $j$ as follows: a tangent vector $V$ at a point $\{p_1,\ldots
,p_r\}\in X_r (f)$ where each $p_i\in f^{-1}(t)$ amounts to a
collection of tangent vectors $v_i\in T_{p_i}X'$ such that all of
the $\pi_* v_i\in T_t S^2$ are the same, and $\mathbb{J}_j V$ is
defined as the collection of vectors $\{jv_1\ldots ,jv_r\}$.  (The
assumption that $f$ is a pseudoholomorphic map with respect to $j$
ensures that the `horizontal parts' $\pi_*jv_i$ all agree, so that
the collection $\{jv_1\ldots ,jv_r\}$ is in fact a well-defined
tangent vector to $X_{r}(f)$; both Section 5 of \cite{Smith} and a
previous version of this paper assert that $\mathbb{J}_j$ can be
constructed if $j$ is merely assumed to make the fibers of $f$
holomorphic, but this is not the case.)  Conversely, a  section
$s$ of $X_r (f)$ naturally corresponds to a closed set $C_s$ in
$X'$ meeting each fiber $r$ times with multiplicities, and $s$ is
$\mathbb{J}_j$--holomorphic exactly if $C_s$ is a $j$--holomorphic
subset of $X'$.  Moreover, as Smith shows, there is just one
homotopy class $c_{\alpha}$ of sections of $X_{r}(f)$ which
tautologically correspond to closed sets in any given class
$\alpha$, and the expected complex dimension $d(\alpha)$ of the
moduli space of such sections is the same as the expected
dimension of the moduli space involved in the construction of the
Gromov invariant. So it seems appropriate to try to count
holomorphic curves in $X$ by counting holomorphic sections of the
various $X_{r}(f)$ in the corresponding homotopy classes.
Accordingly, in \cite{Smith} (and earlier in \cite{DS} for the
special case $\alpha=\kappa$), the \emph{standard surface count}
$\mathcal{DS}_{(X,f)}(\alpha)$ is defined to be the Gromov--Witten
invariant counting sections $s$ of $X_{r}(f)$ in the class
$c_{\alpha}$ with the property that, for a generic choice of
$d(\alpha)$ points $z_i$ in $X$, the value $s(f(z_i))$ is a
divisor in $S^r f(z_i)$ containing the point $z_i$.  Note that
such sections will then descend to closed sets in $X$ containing
each of the points $z_i$. Actually, in order to count curves in
$X$ and not $X'$ $\alpha$ should be a class in $X$, and the
standard surface count will count sections of $X_{r}(f)$ in the
class $c_{\pi_{k}^{*}(\alpha)}$;  it's straightforward to see that
$Gr(\pi_{k}^{*}(\alpha))=Gr(\alpha)$. $k$ here needs to be taken
large enough that the relevant moduli space of sections of
$X_{r}(f)$ is compact; we can ensure that this will be true if $k
[\omega]^2 > \omega \cdot \alpha$, since in this case the section
component of any cusp curve resulting from bubbling would descend
to a possibly-singular symplectic submanifold of $X'$ on which
$\pi_{k}^{*}\omega$ evaluates negatively, which is impossible.
With this compactness result understood, the Gromov--Witten
invariant in question may be defined using the original definition
given by Yongbin Ruan and Gang Tian in \cite{RT}; recourse to
virtual moduli techniques is not necessary.

The main theorem of \cite{Smith}, proven using Serre duality on
the fibers of $f$ and the special structure of the Abel--Jacobi
map from $X_{r}(f)$ to a similarly-defined ``relative Picard
scheme'' $P_{r}(f)$, is that
$\mathcal{DS}_{(X,f)}(\alpha)=\pm\mathcal{DS}_{(X,f)}(\kappa-\alpha)$,
provided that $b^{+}(X)>b_{1}(X)+1$ (and Smith in fact gives at
least a sketch of a proof whenever $b^{+}(X)>2$) and that the
degree of the Lefschetz fibration is sufficiently high.

Smith's theorem would thus provide a new proof of Taubes duality under a somewhat weaker constraint on the Betti numbers if
it were the case that (as Smith conjectures)
\begin{equation} \label{conj} \mathcal{DS}_{(X,f)}(\alpha)=Gr(\alpha) \end{equation}
  Even without this, the duality theorem is strong enough to yield several of the topological consequences of Taubes duality: for
instance, the main theorem of \cite{DS} gives the existence of a
symplectic surface Poincar\'e dual to $\kappa$; see also Section
7.1 of \cite{Smith} for new Seiberg--Witten theory-free proofs of
several other symplectic topological results of the mid-1990s. The
tautological correspondence discussed above would seem to provide
a route to proving the conjecture (\ref{conj}), but one encounters
some difficulties with this. While the tautological correspondence
implies that the moduli space of $\mathbb{J}$--holomorphic
sections of $X_{r}(f)$ agrees set-theoretically with the space of
$j$--holomorphic submanifolds of $X$, it is not obvious whether
the weights assigned to each of the sections and curves in the
definitions of the respective invariants will agree.  This might
seem especially worrisome in light of the fact that the invariant
$Gr$ counts some multiply-covered square-zero tori with weights
other than $\pm 1$ in order to account for the wall crossing that
occurs under a variation of the complex structure when a sequence
of embedded curves converges to a double cover of a square-zero
torus.

This paper confirms, however, that the weights agree.  The main
theorem is:
\begin{thm}\label{main} Let $f\co (X,\omega)\to S^2$ be a symplectic Lefschetz
fibration and $\alpha\in H^2 (X,\mathbb{Z})$ any class such that
$\omega\cdot \alpha <\omega\cdot (fiber)$.  Then
$\mathcal{DS}_{(X,f)}(\alpha)=Gr(\alpha)$. \end{thm}

The hypothesis of the theorem is satisfied, for instance, for
Lefschetz fibrations $f$ of sufficiently high degree obtained by
Donaldson's construction applied to some symplectic manifold $X_0$
($X$ will be a blow-up of $X_0$) where $\alpha$ is the pullback of
some cohomology class of $X_0$.  In particular, the theorem
implies that the standard surface count for such classes is
independent of the degree of the fibration provided that the
degree is high enough. It is not known whether this fact can be
proven by comparing the standard surface counts directly rather
than equating them with the Gromov invariant, though Smith has
suggested that the stabilization procedure discussed in \cite{AK}
and \cite{stab} might provide a route for doing so.

Combining the above Theorem \ref{main} with Theorem \ref{main} of
\cite{Smith}, we thus recover:

\begin{cor}[Taubes] \label{duality}  Let $(X,\omega)$ be a symplectic 4--manifold
with $b^{+}(X)$ $>b_1(X)+1$ and canonical class $\kappa$.  Then for
any $\alpha\in H^{2}(X;\mathbb{Z})$, $Gr(\alpha)=\pm Gr(\kappa
-\alpha)$. \end{cor}

While the requirement on the Betti numbers here is stronger than
that of Taubes (who only needed $b^{+}(X)>1$), the proof of
Corollary \ref{duality} via the path created by Donaldson and
Smith and completed by Theorem \ref{main} avoids the difficult
gauge-theoretic arguments of \cite{book} and also remains more
explicitly within the realm of symplectic geometry.

We now briefly describe the proof of Theorem \ref{main} and the
organization of this paper.  Our basic approach is to try to
arrange to use, for some $j$ making $f$ pseudoholomorphic, the
$j$--moduli space to compute $Gr$ and the $\mathbb{J}_j$--moduli
space to compute $\mathcal{DS}$, and to show that the contribution
of each curve in the former moduli space to $Gr$ is the same as
the contribution of its associated section to $\mathcal{DS}$.  In
Section 2, we justify the use of such $j$ in the computation of
$Gr$.  In Section 3, we refine our choice of $j$ to allow
$\mathbb{J}_j$ to be used to compute $\mathcal{DS}$, at least when
there are no multiple covers in the relevant moduli spaces.  For a
non-multiply-covered curve $C$, then, we show that its
contributions to $Gr$ and $\mathcal{DS}$ agree by, in Section 4,
directly comparing the spectral flows for $C$ and for its
associated section $s_C$ of $X_r (f)$.  This comparison relies on
the construction of an almost complex structure which makes both
$C$ and $f$ holomorphic and which is integrable near $C$. Although
for an arbitrary curve $C$ such an almost complex structure may
not exist, the constructions of Section 3 enable us to reduce to
the case where each curve at issue does admit such an almost
complex structure nearby by first delicately perturbing the
original almost complex structure on $X$. We use this result in
Section 4 to set up corresponding spectral flows in $X$ and $X_r
(f)$ and show that the signs of the spectral flows are the same,
which proves that curves with no multiply-covered components
contribute in the same way to $\mathcal{DS}$ and $Gr$.

For curves with multiply covered components, such a direct
comparison is not possible because the almost complex structure
$\mathbb{J}$ is generally non-differentiable at the image of the
section of $X_r (f)$ associated to such a curve.  Nonetheless, we
see in Section 5 that the contribution of such a $j$--holomorphic
curve $C$ to the invariant $\mathcal{DS}$ is still a well-defined
quantity which remains unchanged under especially nice variations
of $j$ and $C$ and which is the same as the contribution of $C$ to
$Gr$ in the case where $j$ is integrable and nondegenerate in an
appropriate sense. To obtain this contribution, we take a smooth
almost complex structure $J$ which is close in H\"older norm to
$\mathbb{J}$; because Gromov compactness remains true in the
H\"older context, this results in the section $s$ of $X_r (f)$
tautologically corresponding to $C$ being perturbed into some
number (possibly zero) of $J$--holomorphic sections which are
constrained to lie in some small neighborhood of the original
section $s$, and the contribution of $C$ to $\mathcal{DS}$ is then
obtained as the signed count of these nearby sections. We then
deduce the agreement of $\mathcal{DS}$ and $Gr$ by effectively
showing that any rule for assigning contributions of
$j$--holomorphic curves in the 4--manifold $X$ which satisfies the
invariance properties of the contributions to $\mathcal{DS}$ and
agrees with the contributions to $Gr$ in the integrable case must
in fact yield Taubes' Gromov invariant.  Essentially, the fact
that $\mathcal{DS}$ is independent of the almost complex structure
used to define it forces the contributions to $\mathcal{DS}$ to
satisfy wall crossing formulas identical to those introduced by
Taubes for $Gr$ in \cite{Taubes}.  Since the results of Section 3
allow us to assume that our curves admit integrable
 complex structures nearby which make the fibration holomorphic, and we know that
contributions to $\mathcal{DS}$ and $Gr$ are the same in the
integrable case, the wall crossing formulas lead to the result
that $\mathcal{DS}=Gr$ in all cases. This approach could also be
used to show the agreement of $\mathcal{DS}$ and $Gr$ for
non-multiply covered curves, but the direct comparison used in
Section 4 seems to provide a more concrete way of understanding
the correspondence between the two invariants, and most of the
lemmas needed for this direct proof are also necessary for the
indirect proof given in Section 5, so we present both approaches.

Throughout the paper, just as in this introduction, a lowercase
$j$ will denote an almost complex structure on the 4--manifold,
and an uppercase $J$ (or $\mathbb{J}$) will denote an almost
complex structure on the relative Hilbert scheme.  When the
complex structure on the domain of a holomorphic curve appears, it
will be denoted by $i$.

This results of this paper are also contained in my thesis
\cite{thesis}. I would like to thank my advisor Gang Tian for
suggesting this interesting problem and for many helpful
conversations while this work was in progress.

\section{Good almost complex structures I}
Let $f\co X\to S^2$ be a symplectic Lefschetz fibration and
$\alpha\in H^2 (X,\mathbb{Z})$. As mentioned in the introduction,
if $j$ is an almost complex structure on $X$ with respect to which
$f$ is pseudoholomorphic, we have a tautological correspondence
$\mathcal{M}_{X}^{j}(\alpha)=\mathcal{MS}_{X_{r}(f)}^{\mathbb{J}_{j}}(c_{\alpha})$
between the space of $j$--holomorphic submanifolds of $X$
Poincar\'e dual to $\alpha$ with no fiber components and the space
of $\mathbb{J}_j$--holomorphic sections of $X_{r}(f)$ in the
corresponding homotopy class.  In light of this, to show that
$Gr(\alpha)$ agrees with $\mathcal{DS}_{(X,f)}(\alpha)$, we would
like, if possible, to use such an almost complex structure $j$ to
compute the former and the corresponding $\mathbb{J}_j$ to compute
the latter.  Two obstacles exist to carrying this out: first, the
requirement that $j$ make $f$ holomorphic is a rather stringent
one, so it is not immediately clear that the moduli spaces of
$j$--holomorphic submanifolds will be generically well-behaved;
second, the almost complex structure $\mathbb{J}_j$ is only
H\"older continuous, and so does not fit into the general
machinery for constructing Gromov--Witten invariants such as
$\mathcal{DS}$.  The first obstacle will be overcome in this
section.  The second obstacle is more serious, and will receive
its share of attention in due course.

We will, in general, work with Lefschetz fibrations such that
$\omega\cdot\alpha < \omega\cdot (fiber)$ for whatever classes
$\alpha$ we consider; note that this requirement can always be
fulfilled by fibrations obtained by Donaldson's construction, and
ensures that $j$--holomorphic curves in class $\alpha$ never have
any fiber components.

By a \emph{branch point} of a $j$--holomorphic curve $C$ we will
mean a point at which $C$ is tangent to one of the fibers of $f$.

\begin{lemma} \label{good AC} Let $f\co (X,\omega)\to (S^{2},\omega_{FS})$ be a symplectic Lefschetz fibration and let $\alpha\in H^{2}(X,\mathbb{Z})$ be such that
$d=d(\alpha)\geq 0$ and $\omega\cdot\alpha < \omega\cdot (fiber)$.
Let $\mathcal{S}$ denote the set of pairs $(j,\Omega)$ where $j$
is an almost complex structure on $X$ making $f$ holomorphic and
$\Omega$ is a set of $d$ distinct points of $f$, and let
$\mathcal{S}^{0}\subset \mathcal{S}$ denote the set for which:
\begin{enumerate} \item $(j, \Omega)$ is nondegenerate in the sense of Taubes
\cite{Taubes}; in particular, where
$\mathcal{M}_{X}^{j,\Omega}(\alpha)$ denotes the set of
$j$--holomorphic curves Poincar\'e dual to $\alpha$ passing
through all the points of $\Omega$,
$\mathcal{M}_{X}^{j,\Omega}(\alpha)$ is a finite set consisting of
embedded curves. \item Each member of
$\mathcal{M}_{X}^{j,\Omega}(\alpha)$ misses all critical points of
$f$. \item No curve in $\mathcal{M}_{X}^{j,\Omega}(\alpha)$ meets
any of the branch points of any of the other curves.
\end{enumerate} Then $\mathcal{S}^{0}$ is open and dense in
$\mathcal{S}$.\end{lemma}
\begin{proof}
As usual for statements such as the assertion that Condition 1 is
dense, the key is the proof that the map $\mathcal{F}$ defined
from
\[ \mathcal{U}=\{(i,u,j,\Omega)|(j,\Omega)\in \mathcal{S},u\co \Sigma\looparrowright X,  \Omega\subset Im(u), \,u\in W^{k,p}\} \]
to a bundle with fiber $W^{k-1,p}(\Lambda^{0,1}T^* \Sigma \otimes
u^* TX)$ by $(i,u,j,\Omega)\mapsto \bar{\partial}_{i,j}u$ is
submersive at all zeroes. ($i$ denotes the complex structure on
the domain curve $\Sigma$.)

Now as in the proof of Proposition 3.2 of \cite{RT} (but using a
$\dbar$--operator equal to one-half of theirs) , the linearization
at a zero $(i,u,j,\Omega)$ is given by
\[ \mathcal{F}_{*} (\beta,\xi,y,\vec{v})=D_{u}\xi +\frac{1}{2}(y\circ du\circ i + j\circ du\circ \beta) \]
Here $D_u$ is elliptic, $\beta$ is a variation in the complex
structure on $\Sigma$ (and so can be viewed as a member of
$H^{0,1}_i (T_{\mathbb{C}}\Sigma))$ and $y$ is a $j$--antilinear
endomorphism of $TX$ that (in order that $\exp_j y$ have the
compatibility property) preserves $T^{vt}X$ and pushes forward
trivially to $S^2$, so with respect to the splitting
$TX=T^{vt}X\oplus T^{hor}X$ ($T^{hor}$ being the symplectic
complement of $T^{vt}$; of course this splitting only exists away
from $Crit(f)$) $y$ is given in block form as
\begin{displaymath} y=\left( \begin{array}{cc} a & b \\ 0 & 0 \end{array} \right) \end{displaymath}
where all entries are $j$--antilinear.

Now suppose $\eta\in W^{k-1,p}(\Lambda^{0,1}T^* \Sigma \otimes u^*
TX)$, so that $\eta$ is a complex-antilinear map $T\Sigma\to u^*
TX$, and take a point $x_{0}\in \Sigma$ for which $d(f\circ
u)(x_{0})$ is injective.  Let $v$ be a generator for
$T^{1,0}_{x_{0}}\Sigma$; then $du(i(v))\in (T^{1,0}X)_{u(x_{0})}$
and $du(i(\bar{v}))\in (T^{0,1}X)_{u(x_{0})}$ are tangent to
$u(\Sigma)$ and so have nonzero horizontal components.  We take
$y(u(x_{0}))=\left( \begin{array}{cc} 0 & b \\ 0 & 0 \end{array}
\right)$ where
\begin{displaymath} b\co T^{hor}_{u(x_{0})}\to T^{vt}_{u(x_{0})} \end{displaymath}
is a $j$--antilinear map with $b(du(v)^{hor})=(\eta(v))^{vt}$ and
$b(du(\bar{v})^{hor})=(\eta(\bar{v}))^{vt}$.  Since complex
antilinear maps are precisely those maps interchanging $T^{1,0}$
with $T^{0,1}$ this is certainly possible.

Suppose now that $\eta\in \cok (\mathcal{F}_*)_{(i,u,j,\Omega)}$.
The above considerations show that for any point $x_0\notin Crit(f\circ u)$
there is $y$ such that
\begin{equation} \label{Feta} \mathcal{F}_* (0,0,y,0)(x_0)=\eta^{vt}(x_0).\end{equation} Cutting off $y$ by some function
$\chi$ supported near $x_0$, if $\eta^{vt}(x_0)\neq 0$ we can arrange that \[
\int_{\Sigma}\langle \mathcal{F}_* (0,0,\chi y,0),\eta\rangle=\int_{\Sigma}\langle \mathcal{F}_* (0,0,\chi y,0),\eta^{vt}\rangle>0,\]
contradicting the supposition that  $\eta\in \cok (\mathcal{F}_*)_{(i,u,j,\Omega)}$.  $\eta^{vt}$ must
therefore be zero at every point not in $Crit (f\circ u)$.

Meanwhile, letting $\eta^{C}$ denote the projection of $\eta$
(which is an antilinear map $T\Sigma\to u^{*}TX$) to $TC$ where
$C=Im(u)$, $\eta^{C}$ then is an element of the cokernel of the
linearization at $(i,id)$ of the map $(i',v)\mapsto\dbar_{i',i}
v$, $i'$ being a complex structure on $\Sigma$ and $v$ being a map
$\Sigma\to \Sigma$.  But the statement that this cokernel vanishes
is just the statement that the set of complex structures on
$\Sigma$ is unobstructed at $i$ (for the cokernel of the map
$v\to\dbar_{i,i}v$ is $H^{1}(T_{\mathbb{C}}\Sigma)$, which is the
same as the space through which the almost complex structures $i'$
vary infinitesimally, and the relevant linearization just sends a
variation $\beta$ in the complex structure on $\Sigma$ to
$i\beta/2$).  So in fact $\eta^{C}=0$.

Now at any point $x$ on $\Sigma$ at which $(f\circ u)_{*}(x)\neq
0$, $TC$ and $T^{vt}X$ together span $TX$, so since
$\eta^{C}(x)=\eta^{vt}(x)=0$ we have $\eta(x)=0$.  But the
assumption on the size of the fibers ensures that $(f\circ
u)_{*}(x)\neq 0$ for all but finitely many $x$, so $\eta$ vanishes
at all but finitely many $x$, and hence at all $x$ since elliptic
regularity implies that $\eta$ is smooth.  This proves that
$(\mathcal{F}_*)_{(i,u,j,\Omega)}$ is submersive whenever
$\mathcal{F}(i,u,j,\Omega)=0$.  The Sard--Smale theorem applied to
the projection $(i,u,j,\Omega)\mapsto (j,\Omega)$ then gives that
Condition 1 in the lemma is a dense (indeed, generic) condition;
that it is an open condition just follows from the fact that
having excess kernel is a closed condition on the linearizations
of the $\bar{\partial}$, so that degeneracy is a closed condition
on $(j,\Omega)$.

As for Conditions 2 and 3, from the implicit function theorem for
the $\dbar$--equation it immediately follows that both are open
conditions on $(j,\Omega)\in\mathcal{S}$ satisfying Condition 1,
so it suffices to show denseness.  To begin, we need to adjust the
incidence condition set $\Omega$ so that it is disjoint from the
critical locus of $f$ and from all of the branch points of all of
the curves of $\mathcal{M}_{X}^{j,\Omega}(\alpha)$.  So given a
nondegenerate pair $(j,\Omega)$ we first perturb $\Omega$ to be
disjoint from $crit(f)$ while $(j,\Omega)$ remains nondegenerate;
then, supposing a point  $p\in\Omega$ is a branch point of some
$C_0\in\mathcal{M}_{X}^{j,\Omega}(\alpha)$, we change $\Omega$ by
replacing $p$ by some $p'$ on $C_0$ which is not a branch point of
$C_0$ and is close enough to $p$ that for each other curve
$C\in\mathcal{M}_{X}^{j,\Omega}(\alpha)$ which does not have a
branch point at $p$, moving $p$ to $p'$ has the effect of
replacing $C$ in the moduli space by some $C'$ which also does not
have a branch point at $p'$ (this is possible by the implicit
function theorem).  Denoting the new incidence set by $\Omega'$,
the number of curves of $\mathcal{M}_{X}^{j,\Omega'}(\alpha)$
having a branch point at $p'$ is one fewer than the number of
curves of $\mathcal{M}_{X}^{j,\Omega}(\alpha)$  having a branch
point at $p$, and so repeating the process we eventually arrange
that no curve in $\mathcal{M}_{X}^{j,\Omega}(\alpha)$ has a branch
point at any point of $\Omega$.

So now assume $(j,\Omega)\in\mathcal{S}$ with $\Omega$ missing
both $Crit(f)$ and all branch points of all curves in
$\mathcal{M}_{X}^{j,\Omega}(\alpha)$.  Let \[
\mathcal{M}_{X}^{j,\Omega}(\alpha)=\{[u_{1}],\ldots,[u_{r}]\} \]
where $[u_m]$ denotes the equivalence class of a map $u_m$ under
the action of $Aut(\Sigma_m)$, $\Sigma_{m}$ being the (not
necessarily connected) domain of $u_{m}$.  For each $m$, enumerate
the points of $\Sigma_m$ which are mapped by $u_m$ either to $Crit
(f)$ or to an intersection point with one of the other curves as
$p_{m,1},\ldots,p_{m,l}$, so in particular none of the $u_m
(p_{m,k})$ lie in $\Omega$. Take small, disjoint neighborhoods
$U_{m,k}$ of the $p_{m,k}$ such that $u_{m}(U_{m,k})$ misses
$\Omega$ and $u_{m}(U_{m,k}\setminus\frac{1}{2}U_{m,k})$ misses
each of the other curves and also misses $Crit(f)$, and take
\emph{local} sections $\xi_{m,k}$ of $u_{m}^{*}T^{vt}X$ over
$U_{m,k}$ such that $D_{u_{m}}\xi_{m,k}=0$ and
$\xi_{m,k}(p_{m,k})\neq 0$ (this is certainly possible, as the
$\xi_{m,k}$ only need to be defined on small discs, on which the
equation $D_{u_{m}}\xi_{m,k}=0$ has many solutions). Now for each
$m$ glue the $\xi_{m,k}$ together to form $\xi_{m}\in
\Gamma(u_{m}^{*}T^{vt}X)$ by using cutoff functions which are 1 on
$\frac{1}{2}U_{m,k}$ and 0 outside $U_{m,k}$.  Then since
$D_{u_{m}}\xi_{m,k}=0$ the sections $D_{u_{m}}\xi_{m}$ will be
supported in
\[ A_{m}=\bigcup_{k}(U_{m,k}\setminus \frac{1}{2}U_{m,k}). \]
Now according to page 28 of \cite{MS}, the linearization $D_{u_m}$
may be expressed with respect to a $j$--Hermitian connection
$\nabla$ by the formula \begin{equation} \label{Du}
(D_{u_m}\xi)(v)=\frac{1}{2}(\nabla_v \xi +
j(u_m)\nabla_{iv}\xi)+\frac{1}{8}N_j((u_{m})_{*}v,\xi)
\end{equation} where $N_j$ is the Nijenhuis tensor.  Our sections
$\xi_m$ are vertically-valued, so the first two terms above will
be vertical tangent vectors; in fact, the last term will be as
well, because where $z$ is the pullback of the local coordinate on
$S^2$ and $w$ a holomorphic coordinate on the fibers, the
anti-holomorphic tangent space for $j$ can be written \[ T^{0,1}_j
X=\langle\partial_{\bar{z}}+b(z,w)\partial_{w},\partial_{\bar{w}}
\rangle, \] in terms of which one finds \begin{equation}
\label{Nj}
N_{j}(\partial_{\bar{z}},\partial_{\bar{w}})=4(\partial_{\bar{w}}b)\partial_w.
\end{equation}  So if $\xi$ is a vertically-valued vector field, the right-hand side of Equation\break \ref{Du}
is also vertically-valued for any $v$, ie, $D_{u_m}$ maps
$W^{k,p}(u_{m}^{*}T^{vt}X)$ to\break
$W^{k-1,p}(\Lambda^{0,1}T^{*}\Sigma_m\otimes u_{m}^{*}T^{vt}X)$
(and not just to $W^{k-1,p}(\Lambda^{0,1}T^{*}\Sigma_m\otimes
u_{m}^{*}TX)$). Now \[ D_{u_m}\xi_{m}\in
W^{k-1,p}(\Lambda^{0,1}T^{*}\Sigma_m\otimes u_{m}^{*}T^{vt}X) \]
is supported in $A_m$, so (using that $u_m(A_m)$ misses $Crit
(f)$) as in (\ref{Feta}) we can find a perturbation $y_m$ of the
almost complex structure $j$ supported near $u_{m}(A_m)$ such that
\[
\mathcal{F}_{*}(0,\xi_{m},y_{m},0)=D_{u_{m}}\xi_{m}+\frac{1}{2}y_{m}\circ
du_{m}\circ m=0.\]
Since the $u_{m}(\bar{A_{m}})$ are disjoint, we can paste these
$y_{m}$ together to obtain a global perturbation $y$ with
$\mathcal{F}_{*}(0,\xi_{m},y,0)=0$ for each $m$.  For $t>0$ small
enough that $(\exp_{j}(ty),\Omega)$ remains nondegenerate, the
holomorphic curves for the complex structure $\exp_{j}(ty)$ will
be approximated in any $W^{k,p}$ norm ($p>2$) to order
$C\|\exp_{j}(ty)-j\|_{C^{1}}\|t\xi_{m}\|_{W^{k,p}}\leq Ct^{2}$ by
the curves $\exp_{u_{m}}(t\xi_{m})$ (using, for example, the
implicit function theorem as formulated in Theorem 3.3.4 and
Proposition 3.3.5 of \cite{MS}).  Now since $\xi_m (p_{m,k})\neq
0$, the $\exp_{u_{m}}(t\xi_{m})$ will have their branch points
moved vertically with respect to where they were before; in
particular, these curves will no longer pass through $Crit(f)$,
and their branch points will no longer meet other curves.
Similarly (for $t$ suitably small, and $k$ appropriately large
chosen at the beginning of the procedure) any set of curves within
$Ct^{2}$ of these in $W^{k,p}$--norm will satisfy these conditions
as well.  So for $t$ small enough, $(\exp_{j}(ty),\Omega)$ will
obey conditions 1 through 3 of the lemma. $(j,\Omega)$ was an
arbitrary nondegenerate pair, so it follows that $\mathcal{S}^{0}$
is dense.  \end{proof}

As has been mentioned above, the almost complex structure
$\mathbb{J}_j$ that we would in principle like to use to evaluate
$\mathcal{DS}$ is generally only H\"older continuous; however,
under certain favorable circumstances we shall see presently that
it is somewhat better-behaved. To wit, assume that our almost
complex structure $j$ is given locally by \[
T^{0,1}_{j}=\langle\partial_{\bar{z}}+b(z,w)\partial_{w},\partial_{\bar{w}}
\rangle, \] where $z$ is the pullback of the coordinate on the
base and $w$ a coordinate on the fibers.  Then, following
\cite{ST2}, where $\sigma_k$ denotes the $k$th elementary
symmetric polynomial, the function
\[
\hat{b}_{d}(z,w_1,\ldots,w_r)=\sum_{k=1}^{r}\sigma_{d-1}(w_1,\ldots,\widehat{w_k},\ldots,w_r)b(z,w_k)
\] on $\mathbb{C}\times \mathbb{C}^r$ is symmetric in the $w_k$
and so descends to a function $b_d (z,\sigma_1,\ldots,\sigma_r)$
on $\mathbb{C}\times S^r\mathbb{C}$, and our almost complex
structure $\mathbb{J}_j$ on $X_r (f)$ is given locally by \[
T^{0,1}_{\mathbb{J}_j}=\langle\partial_{\bar{z}}+\sum_{d=1}^{r}b_{d}(z,\sigma_1,\ldots,\sigma_r)\partial_{\sigma_d},\partial_{\bar{\sigma}_1},\ldots,
\partial_{\bar{\sigma}_r} \rangle.\]
The nondifferentiability of $\mathbb{J}_j$ can then be understood
in terms of the fact that smooth symmetric functions on
$\mathbb{C}^r$ such as $\hat{b}_d (z,\cdot)$ generally only
descend to H\"older continuous functions in the standard
coordinates $\sigma_1,\ldots,\sigma_r$ on $S^{r}\mathbb{C}$ (when
$r=2$, for example, consider the function
$\bar{w}_1w_2+w_1\bar{w}_2$).  On the other hand,
\emph{holomorphic} symmetric functions on $\mathbb{C}^r$ descend
to holomorphic (and in particular smooth) functions on the
symmetric product, so when $\partial_{\bar{w}}b=0$, the functions
$b_d$ are holomorphic in the vertical coordinates, and so
$\mathbb{J}_j$ is smooth. Furthermore, note that by Equation
\ref{Nj}, $b$ is holomorphic in $w$ exactly when $j$ is integrable
on the neighborhood under consideration; moreover, computing the
Nijenhuis tensor of $\mathbb{J}_j$ shows that $\mathbb{J}_j$ is
integrable exactly when $\partial_{\bar{\sigma}_k}b_l=0$ for all
$k$ and $l$.  This sets the stage for the following proposition,
which foreshadows some of the constructions in the next two
sections:
\begin{prop}\label{Jjint} Let $C\in \mathcal{M}_{X}^{j,\Omega}(\alpha)$ where
$(j,\Omega)$ is as in Lemma \ref{good AC}, and let $s_C$ be the
corresponding section of $X_r (f)$.  If $j$ is integrable on a
neighborhood of $C$, then $\mathbb{J}_j$ is integrable on a
neighborhood of $s_C$.  More generally, if $j$ is only integrable
on neighborhoods of each of the branch points of $C$, then
$\mathbb{J}_j$ is still smooth on a neighborhood of $s_C$.\end{prop}
\begin{proof}
The first statement follows directly from the above argument. As
for the second statement, note that the only place where our
functions $b_d$ above ever fail to be smooth is in the diagonal
stratum $\Delta$ of $\mathbb{C}\times S^{r}\mathbb{C}$ where two
or more points in the divisor in $S^r\mathbb{C}$ come together. A
suitably small neighborhood of $s_C$ only approaches this stratum
in a region whose differentiable structure for the vertical
coordinates is just that of the Cartesian product of symmetric
products of neighborhoods of all the branch points in some fiber
(where smoothness is taken care of by the integrability
assumption) with copies of $\mathbb{C}$ corresponding to
neighborhoods of each of the other points of $C$ which lie in the
same fiber.
\end{proof}

We close this section with a proposition which shows that if
$\mathbb{J}_j$ can be assumed smooth, then its moduli spaces will
generically be well-behaved. We make here a statement about
generic almost complex structures from a set $\mathcal{S}^1$ which
at this point in the paper has not yet been proved to be nonempty;
rest assured that it will be seen to be nonempty in the following
section.

\begin{prop} \label{JjOK} For generic $(j,\Omega)$ in the set $\mathcal{S}^1$ consisting of members of the set $\mathcal{S}^0$
from Lemma \ref{good AC} which satisfy the additional property
that $j$ is integrable near every branch point of every curve $C$
in $\mathcal{M}_{X}^{j,\Omega}(\alpha)$, the linearization of the
operator $\dbar_{\mathbb{J}_j}$ is surjective at each of the
sections $s_C$.  \end{prop}
\begin{proof}
We would like to adapt the usual method of constructing a
universal moduli space
$\mathcal{U}=\{(s,j,\Omega)|\bar{\partial}_{\mathbb{J}_{j}}u=0,
(j,\Omega)\in \mathcal{S}^1,\Omega\subset C_s\}$, appealing to the
implicit function theorem to show that $\mathcal{U}$ is a smooth
Banach manifold, and then applying the Sard--Smale theorem to the
projection from $\mathcal{U}$ onto the second factor (ie,
$\mathcal{S}^1$) to obtain the statement of the proposition.  Just
as in the proof of Lemma \ref{good AC}, this line of argument will
work as long as we can show that the map $(s, j,\Omega)\mapsto
\bar{\partial}_{\mathbb{J}_{j}} s$ is transverse to zero.

Arguing as before, it's enough to show that, for a section $s$
with $\bar{\partial}_{\mathbb{J}_{j}} s=0$, where $D_{s}^{*}$
denotes the formal adjoint of $D_{s}$, and where $i$ denotes the
complex structure on $S^2$, the following holds: if
$D_{s}^{*}\eta=0$, and if, for every variation $y$ in the complex
structure $j$ on $X$ among almost complex structures $j'$ with
$(j',\Omega)\in \mathcal{S}^1$, letting $Y$ denote the variation
in $\mathbb{J}_{j}$ induced by $y$, we have that
\begin{equation} \label{cxvar}
 \int_{S^2} \langle \eta, Y(s)\circ ds \circ i \rangle =0, \end{equation}
then $\eta\equiv 0$.   If $\eta$ were nonzero, then it would be
nonzero at some $t_0\in S^2$ which is not the image under $f$ of
any of the branch points of $C_s$, so assume this to be the case.
Now $\eta$ is a $s^{*}T^{vt}X_{r}(f)$--valued (0,1)--form, so
giving its value at $t_{0}$ is equivalent to giving $r$ maps
$\eta_k\co T_{t_0}S^2\to T_{s_k(t_0)}^{vt}X$ $(r=1,\ldots ,k)$,
where the $s_{k}(t_{0})$ are the points in the fiber
$\Sigma_{t_{0}}$ over $t_{0}$ of the Lefschetz fibration which
correspond to the point $s(t_{0})\in S^{r}\Sigma_{t_{0}}$ (our
assumption on $t_0$ ensures that these are all distinct).
$\eta(t_{0})$ being nonzero implies that one of these cotangent
vectors (say $\eta_{m}$) is nonzero. Then $s_m$ is a local
holomorphic section of $X\to S^2$ around $t_0$, and exactly as in
the proof of Lemma \ref{good AC} we may find a perturbation $y_0$
of the almost complex structure near $s_m (t_0)$ such that \[ y_0
(s_m (t_0))\circ ds_m (t_0)\circ i = \eta_m \] and $y_0$ preserves
the pseudoholomorphicity of the fibration $f$.  Multiplying $y_0$
by a smooth cutoff supported in a suitably small neighborhood of
$s_{m}(t_{0})\in X$, we obtain a variation $y$ of the complex
structure on $X$ whose associated variation $Y$ in
$\mathbb{J}_{j}$ violates (\ref{cxvar}); note that since $y$ is
supported away from the nodes of the curves of
$\mathcal{M}_{X}^{j,\Omega}(\alpha)$, the variation will also not
disrupt the integrability condition in the definition of
$\mathcal{S}^1$. This contradiction shows that $\eta$ must vanish
everywhere, and hence that  $(s, j,\Omega)\mapsto
\bar{\partial}_{\mathbb{J}_{j}} s$ is indeed transverse to zero,
so that the universal space $\mathcal{U}$ will be a manifold and
the usual Sard--Smale theorem argument implies the proposition.
\end{proof}

\section{Good almost complex structures II} \label{gac2}
We fix a symplectic Lefschetz fibration $f\co X\to S^{2}$ and a
class $\alpha\in H^{2}(X,\mathbb{Z})$.  Assume unless otherwise
stated that $(j,\Omega)\in \mathcal{S}^{0}$, so that each curve
$C\in\mathcal{M}_{X}^{j,\Omega}(\alpha)$ is identified by the
tautological correspondence with a section $s_{C}$ of $X_{r}(f)$
which misses the critical locus.  Assume also that $\alpha$ cannot
be decomposed as a sum of classes each of which pairs positively
with $\omega$ and one of which, say $\beta$, satisfies
$\kappa\cdot\beta=\beta\cdot\beta=0$. Then the contribution of
$C\in\mathcal{M}_{X}^{j,\Omega}(\alpha)$ to the invariant
$Gr(\alpha)$ is found by considering a path of operators $D_{t}$
acting on sections of the disc normal bundle $U_{C}$ of $C$ such
that $D_{0}$ is the $\bar{\partial}$ operator obtained from the
complex structure $j_{0}$ on $U_{C}$ given by pulling back
$j|_{C}$ to $U_{C}$ via the Levi--Civita connection, while $D_{1}$
is the $\bar{\partial}$ operator obtained by viewing $U_{C}$ as a
tubular neighborhood of $C$ in $X$ and restricting $j$ to $U_{C}$
(see section 2 of \cite{Taubes}).  If the path $(D_{t})$ misses
the stratum of operators with 2-dimensional kernel and meets the
stratum with one-dimensional kernel transversely, then the
contribution of $C$ to $Gr(\alpha)$ is given by $-1$ raised to a
power equal to the number of times it meets this latter stratum;
more generally the contribution is found by orienting the
zero-dimensional space $\ker D_{1}$ so that the corresponding
orientation of $\det(D_{1})=\Lambda^{max}\ker
D_{1}=\Lambda^{max}\ker D_{1}\otimes(\Lambda^{max} \cok
D_{1})^{*}$ agrees with the natural orientation of the bundle
$\bigcup_{t} \det(D_{t})\times \{t\}$ which restricts to $t=0$ as
the complex orientation of $\det(D_{0})$ (since $j_{0}$ is
integrable, one has
\begin{equation} \label{D} D_{0}\xi=\frac{1}{2}(\nabla \xi + j(u) \nabla \xi \circ i) + \frac{1}{8}N_{j_{0}}(\partial_{j}u,\xi)=\frac{1}{2}(\nabla \xi + j(u) \nabla \xi \circ i) \end{equation}
where $u\co (\Sigma,i)\to X$ is an embedding of $C$, $\nabla$ is a
$j$--hermitian connection, and $N$ is the Nijenhuis tensor, using
remark 3.3.1 of \cite{MS}.  $D_{0}$ therefore commutes with
$j_{0}$, giving $\det(D_{0})$ a natural (complex) orientation).

As for $\mathcal{DS}$, if $J$ is a smooth regular almost complex
structure on $X_{r}(f)$ and
$s\in\mathcal{MS}_{X_r(f)}^{J,\Omega}(c_{\alpha})$, the
contribution of $s$ to $\mathcal{DS}_{(X,f)}(\alpha)$ is similarly
obtained by the spectral flow. Owing to the tautological
correspondence, we would prefer to replace this smooth $J$ with
the almost complex structure $\mathbb{J}_j$.  In general this is
problematic because of the nondifferentiability of $\mathbb{J}_j$,
but let us suppose for a moment that we have found some way to get
around this, by choosing $j$ as in Proposition \ref{JjOK}.
$\mathbb{J}_j$ is then smooth and nondegenerate (ie, the
linearization of $\dbar_{\mathbb{J}_j}$ is surjective) at each of
the sections in the set
$\mathcal{MS}_{X_{r}(f)}^{\mathbb{J}_{j},\Omega}(c_{\alpha})$ of
$\mathbb{J}_j$--holomorphic sections descending to curves which
pass through $\Omega$, which makes the following simple
observation relevant.

\begin{prop}\label{obvious} Assume $J$ is an almost complex structure on
${X}_r (f)$ which is H\"older continuous globally and
smooth and nondegenerate at each member $s$ of
$\mathcal{MS}_{X_{r}(f)}^{J,\Omega}(c_{\alpha})$.  Then
$\mathcal{DS}_{(X,f)}(\alpha)$ may be computed as the sum of the
spectral flows of the linearizations of $\dbar_J$ at the sections
$s$. \end{prop}
\begin{proof} If $J$ were globally smooth this would just be the
definition of $\mathcal{DS}$.  As it stands, we can find a
sequence of smooth almost complex structures $J_n$ agreeing with
$J$ on an open subset $U$ of its smooth locus which contains the
images of all members of
$\mathcal{MS}_{X_{r}(f)}^{J,\Omega}(c_{\alpha})$ such that $J_n$
converges to $J$ in H\"older norm. According to \cite{Sikorav},
Gromov compactness holds assuming only H\"older convergence of the
almost complex structures, so since there are no sections in
$\mathcal{MS}_{X_{r}(f)}^{J,\Omega}(c_{\alpha})$ meeting
$X_{r}(f)\setminus \bar{U}$, for large enough $n$ there must not
be any sections in
$\mathcal{MS}_{X_{r}(f)}^{J_n,\Omega}(c_{\alpha})$ meeting that
region either.  But then since $J_n$ agrees with $J$ on $U$, we
must have $\mathcal{MS}_{X_{r}(f)}^{J_n,\Omega}(c_{\alpha})=
\mathcal{MS}_{X_{r}(f)}^{J,\Omega}(c_{\alpha})$.  Moreover, the
spectral flow for a $J'$--holomorphic section $s$ depends only on
the restriction of $J'$ to a neighborhood of $s$, so since $J$ and
$J_n$ agree near all members of
$\mathcal{MS}_{X_{r}(f)}^{J,\Omega}(c_{\alpha})$, they will both
give the same spectral flows. Using $J_n$ to compute
$\mathcal{DS}$ then proves the proposition. \end{proof}

Assuming then that we can contrive to use the almost complex
structure $\mathbb{J}_j$ to compute $\mathcal{DS}$, we would like
to arrange that the spectral flows for $j$ on the disc normal
bundle and for $\mathbb{J}_j$ on the disc bundle in
$s_{C}^{*}T^{vt}X_{r}(f)$ correspond in some natural way. Now
since $D_{0}$ on $U_{C}\subset X$ comes from a complex structure
which does not preserve the fibers of $f$ (rather, it preserves
the fibers of the normal bundle) and so does not naturally
correspond to any complex structure on a neighborhood of
$Im(s_{C})$ in $X_{r}(f)$, this at first seems a tall order.
However, the key observation is that rather than starting the
spectral flow at $D_{0}$ we can instead start it at the
$\bar{\partial}$ operator $\tilde{D}$ corresponding to any
integrable complex structure $\tilde{j}$ on $U_{C}$. Indeed, if
$j_{t}$ is a path of (not-necessarily integrable ) almost complex
structures from $j_{0}$ to $\tilde{j}$ then the operators
$D_{t}\xi=\frac{1}{2}(\nabla^{t} \xi + j(u) \nabla^{t} \xi \circ
i)$ ($\nabla^{t}$ being a $j_{t}$--Hermitian
 connection) form a family of complex linear operators which by
(\ref{D})  agree at the endpoints with $D_{0}$ and $\tilde{D}$, so
the complex orientation of $\bigcup \det(D_{t})\times \{t\}$
agrees at the endpoints of $D_{0}$ and $\tilde{D}$.  So by taking
the path used to find the contribution of $C$ to $Gr$ to have
$D_{1/2}=\tilde{D}$, the orientation induced on $\det(D_{1})$ by
$\bigcup_{t\in [0,1]} \det(D_{t})\times \{t\}$ and the complex
orientation on $\det(D_{0})$ is the same as that induced by
$\bigcup_{t\in [1/2,1]} \det(D_{t})\times \{t\}$ and the complex
orientation of $\det(D_{1/2})=\det(\tilde{D})$.

The upshot is that for both $Gr$ and $\mathcal{DS}$ we can obtain
the contribution of a given curve (or section) by starting the
spectral flow at any complex structure which is integrable on a
neighborhood of the curve (or section) and makes the curve (or
section) holomorphic.  By Proposition \ref{Jjint}, if $\tilde{j}$
makes $f$ pseudoholomorphic and is integrable on an open set
$U\subset X$ then the corresponding almost complex structure
$\mathbb{J}_{\tilde{j}}$ is integrable on the corresponding
neighborhood in $X_{r}(f)$.  So if we can take $(j,\Omega)$ to
belong to the set $\mathcal{S}^1$ of Proposition \ref{JjOK} (a set
we have not yet shown to be nonempty), we can hope to have the
spectral flows correspond if we can find an almost complex
structure $\tilde{j}$ integrable on a neighborhood of any given
member $C$ of $\mathcal{M}_{X}^{j,\Omega}(\alpha)$ which makes
both $C$ and $f$ holomorphic.  We will see later on that given
such a $(j,\Omega)\in \mathcal{S}^1$, constructing $\tilde{j}$ is
fairly easy, so we turn now to the task of replacing our original
pair $(j,\Omega)$, assumed to be as in Lemma \ref{good AC}, by a
pair belonging to $\mathcal{S}^1$.

Accordingly, let $C\in \mathcal{M}_{X}^{j,\Omega}(\alpha)$ where
$(j,\Omega)\in \mathcal{S}^{0}$, and let $u\co \Sigma\to X$ be an
embedding of $C$.  Restrict attention to a small neighborhood $U$
of a branch point $p$ of $C$; note that by Condition 3 of Lemma
\ref{good AC}, $U$ may be taken small enough to miss all of the
other curves in $\mathcal{M}_{X}^{j,\Omega}(\alpha)$; also, as is
shown in the proof of that Lemma, $U$ can be taken small enough to
miss $\Omega$. Let $w$ be a $j$--holomorphic coordinate on the
fibers, and let $z$ be the pullback of the holomorphic coordinate
on the base $S^{2}$, translated so that $p$ has coordinates
$(0,0)$. Then $j$ is determined by giving a function $b$ such that
the anti-holomorphic tangent space for $j$ is
\begin{equation} \label{T01} T^{0,1}_{j}=\langle \partial_{\bar{z}}+b(z,w)\partial_{w},\partial_{\bar{w}}\rangle  \end{equation}
From Equation \ref{Nj}, a complex structure defined by such an
expression is integrable exactly when $b_{\bar{w}}\equiv 0$.

In general, we cannot hope to realize our initial goal of finding
an almost complex structure making both $f$ and $C$ holomorphic
which is integrable on a neighborhood of $C$. The problem may be
explained as follows.  If our almost complex structure is to have
the form (\ref{T01}), the condition that $C$ be holomorphic
determines $b|_{C}$ uniquely. In regions not containing any points
of $crit(f|_{C})$ this doesn't create a problem, since at least
after shrinking the region so that each connected component of its
intersection with any fiber contains only one point of $C$,
$b|_{C}$ can be extended to the region arbitrarily, say by
prescribing $b$ to be locally constant on each fiber.  When $C$ is
tangent to the fiber $\{w=0\}$ at $(0,0)$, though, we have that
$\partial_{\bar{w}}\in T_{(0,0)}C\otimes\mathbb{C}$, and so
$b_{\bar{w}}(0,0)$ is determined by $b|_{C}$ (which is in turn
determined by $C$).

More concretely, assuming the tangency between $C$ and the fiber
at $(0,0)$ to be of second order, we can write $C=\{z=g(w)\}$
where, after scaling $w$, $g$ is a function of form
$g(w)=w^{2}+O(3)$.  A routine computation shows that for $C$ to be
holomorphic with respect to an almost complex structure defined by
(\ref{T01}), we must have
\begin{equation} \label{b} b(g(w),w)=\frac{-g_{\bar{w}}}{|g_{w}|^{2}-|g_{\bar{w}}|^{2}} \end{equation}
from which one finds by using the Taylor expansion of $g$ to
Taylor-expand the right-hand side that
$b_{\bar{w}}(0,0)=-\frac{1}{8}g_{w\bar{w}\bar{w}\bar{w}}(0)$,
which has no a priori reason to be zero.

Evidently, then, in order to construct an almost complex structure
$\tilde{j}$ as above, or even to find a pair $(j_1,\Omega)\in
\mathcal{S}^1$, so that $j_{1}$ is integrable in neighborhoods of
all of the branch points of all of the curves in
$\mathcal{M}^{j_{1},\Omega}_{X}(\alpha)$, we will have to move the
$j$--holomorphic curves $C$. We show now how to arrange to do so.

Let $j$, $\Omega$, $C$, $u$, $p$, and $U$ be as above.  We will
construct almost complex structures $j_{\epsilon}$ which are
integrable on increasingly small neighborhoods of $p$ and the
linearization of whose $\bar{\partial}$ operators (considered as
acting on sections of the normal bundle $N=N_{C}=N_{C}X$) are
increasingly close to the linearization of $\bar{\partial}_{j}$.
For the latter condition one might initially expect that the
$j_{\epsilon}$ would need to be $C^{1}$--close to $j$, which the
above considerations indicate would be impossible in the
all-too-likely event that $b_{\bar{w}}(0,0)\neq 0$.  However, the
only directional derivatives of the complex structure which enter
into the formula for the linearization are those in the direction
of the section being acted on, so since normal vectors of $C$ near
$p$ have small vertical components the disagreement between the
vertical derivatives of $j_{\epsilon}$ and $j$ will turn out not
to pose a problem.

To begin, we fix $r$ and $\epsilon_{0}$ such that the set \[
D_{3r}^z \times D_{3\epsilon_{0}}^w :=\{(z,w)\mid |z|<3r,
|w|<3\epsilon_{0}\} \] is disjoint from all curves of
$\mathcal{M}_{X}^{j,\Omega}(\alpha)$ except for $C$.  Let $\beta
(z)$ (resp. $\chi (w)$) be a cutoff function which is 1 on
$D_{r}^z$ (resp. $D_{\epsilon_{0}}^w$) and 0 outside $D_{2r}^{z}$
(resp.  $D_{2\epsilon_{0}}^w$).  Let
\[ C_{0}=\sup\{|\nabla\beta|,|\nabla\chi|/\epsilon_{0}\} \] (so we can certainly take $C_{0}\leq\max\{2/r,2\}$).  Where
\[T^{0,1}_{j}=\langle\partial_{\bar{z}}+b(z,w)\partial_{w},\partial_{\bar{w}} \rangle\] for each $\epsilon<\epsilon_{0}$ we define almost complex structures $j_{\epsilon}$ by
\begin{equation} \label{jep}  T^{0,1}_{j_{\epsilon}}=\langle \partial_{\bar{z}}+b_{\epsilon}(z,w)\partial_{w},\partial_{\bar{w}}\rangle
\end{equation}
where
\[  b_{\epsilon}(z,w)=\beta(z)\chi\left(\frac{\epsilon_{0}w}{\epsilon}\right)b(z,0)+\left(1-\beta(z)\chi\left(\frac{\epsilon_{0}w}{\epsilon}\right)\right)b(z,w) \]
So within the region $D_{r}^{z}\times D_{\epsilon}^{w}$ we have
$(b_{\ep})_{\bar w}\equiv 0$, meaning that $j_{\epsilon}$ is
integrable, while outside the region $D_{2r}^{z}\times
D_{2\epsilon}^{w}$ $j_{\epsilon}$ agrees with $j$.  Further,
\begin{equation} |b(z,w)-b_{\epsilon}(z,w)|=|\beta (z)\chi (\epsilon_{0}w/\epsilon)(b(z,w)-b(z,0))|\leq 2\epsilon\|b\|_{C^{1}} \end{equation} (since the expression is zero for $|w|>2\epsilon$),
\begin{align} \label{delz} |\nabla_{z}(b-b_{\epsilon})|&\leq |\nabla_{z}\beta||\chi(\epsilon_{0}w/\epsilon)(b(z,w)-b(z,0))|\notag\\
&\hspace{5cm}+\beta(z)\chi(\epsilon_{0}w/\epsilon)|\nabla_{z}(b(z,w)-b(z,0)| \notag \\
&\leq 2C_{0}\epsilon\|b\|_{C^{1}} + 2\epsilon\|b\|_{C^{2}}
\end{align}
and \begin{align} \label{delw} |\nabla_{w}(b-b_{\epsilon})|&\leq|\nabla_{w}\chi (\epsilon_{0}w/\epsilon)||\beta (z)(b(z,w)-b(z,0))|+\beta(z)\chi(\epsilon_{0} w/\epsilon)|\nabla_{w}b(z,w)| \notag \\
&\leq \frac{C_{0}}{\epsilon}2\epsilon\|b\|_{C^{1}} +
\|b\|_{C^{1}}=(2C_{0}+1)\|b\|_{C^{1}} \end{align}
$C$ is tangent to $\{w=0\}$ at $(0,0)$, so after scaling $z$ we
can write $C$ as $\{z=w^{n}+O(n+1)\}$ for some $n>1$.  It follows
that there is a constant $C_{1}$ such that if $\xi$ is a normal
vector to $C$ based at $(z,w)\in C$ then $|\xi^{vt}|\leq
C_{1}|w|^{n-1}|\xi|$. Hence since
$|\partial_{\xi}(b_{\epsilon}-b)|=0$ if $|w|>2\epsilon$, equations
(\ref{delz}) and (\ref{delw}) give that
\begin{align} \label{delxi} |\partial_{\xi}(b_{\epsilon}-b)| & \leq |\xi^{hor}||\nabla_{z}(b_{\epsilon}-b)|+|\xi^{vt}||\nabla_{w}(b_{\epsilon}-b)| \notag \\
&\leq 2(C_{0}\|b\|_{C^{1}}+\|b\|_{C^{2}})\epsilon |\xi| +
(2C_{0}+1)\|b\|_{C^{1}}C_{1}(2\ep)^{n-1}|\xi| \end{align} We
summarize what we have found in:
\begin{lemma}\label{jeplemma} The almost complex structures given by (\ref{jep}) are integrable in $D_{r}^{z}\times D_{\epsilon}^{w}$ and agree with $j$ outside $D_{2r}^{z}\times D_{2\epsilon}^{w}$.  Further
there is a constant $C_{2}$ depending only on $j$ and the curve $C$ such that $\|j_{\epsilon}-j\|_{C^{0}}\leq C_{2}\epsilon$ and, for any vector $\xi$ normal to $C$,
$|\partial_{\xi}j_{\epsilon}-\partial_{\xi}j|\leq C_{2}\epsilon |\xi|$. \end{lemma}

Now for any almost complex structure $J$ on $X$, the linearization
of $\bar{\partial}_{J}$ at a map \\ $u\co (\Sigma,i)\to (X,J)$ is
given by
\[ D_{u}^{J}\xi=\frac{1}{2}(\nabla^{J}\xi +J(u)\circ \nabla^{J}\xi\circ i)+\frac{1}{2}(\nabla^{J}_{\xi}J)(u)\partial_{J}(u)\circ i \]
where $\nabla^{J}$ is the Levi--Civita connection of the metric
associated to $J$ (this is equation 3.2 of \cite{MS}; they view
$D$ as acting on sections of $u^{*}TX$, but we may equally well
view it as a map $\Gamma(u^{*}N_C)\to\Gamma(u^{*}N_C \otimes
T^{0,1}C)$, as in \cite{Taubes}).  Now the difference between
$\nabla^{j_{\epsilon}}$ and $\nabla^{j}$ is controlled by the
$C^{0}$ norm of $j_{\epsilon}-j$, as is
$\partial_{j_{\epsilon}}(u)-\partial_{j}(u)$, so in the only terms
in which the derivatives of $j_{\epsilon}$ and $j$ come into play
in $(D^{j_{\epsilon}}_{u}-D^{j}_{u})\xi$, the complex structure is
being differentiated in the direction $\xi$.  Lemma \ref{jeplemma}
thus implies:
\begin{cor} \label{lins} There is a constant $C_{3}$ such that the linearizations
\[ D^{j_{\epsilon}}_{u},D^{j}_{u}\co W^{1,p}(u^{*}N_{C})\to L^{p}(u^{*}N_{C}\otimes T^{0,1}C) \]
obey $\|D^{j_{\epsilon}}_{u}\xi-D^{j}_{u}\xi\|_{L^{p}}\leq C_{3}\epsilon \|\xi\|_{W^{1,p}}$. \end{cor}

Now let $D^{\epsilon}$ denote the operator \[
D^{j_{\epsilon}}_{u}\oplus(ev_{\Omega})_{*}\co
W^{1,p}(u^{*}N_{C})\to L^{p}(u^{*}N_{C}\otimes
T^{0,1}C)\oplus\bigoplus_{q\in \Omega}T_{q}X \] and likewise
$D=D^{j}_{u}\oplus (ev_{\Omega})_{*}$.  $D$ and all of the
$D^{\epsilon}$ are then Fredholm of index zero, and $j$ being
nondegenerate in the sense of Taubes \cite{Taubes} amounts to the
statement that $D$ is surjective and hence has a two-sided (since
$\ind (D)=0$) bounded inverse, which we denote $Q$.

\begin{lemma} \label{xin}  Let $\epsilon_{n}\to 0$ and let $\xi_{n}$ be a bounded sequence in $W^{1,p}(u^{*}N_{C})$ with $D^{\epsilon_{n}}\xi_{n}\to 0$.  Then $\xi_{n}\to 0$. \end{lemma}
\begin{proof} The proof is based on the elliptic estimate \begin{equation} \label{ell} \|\xi\|_{W^{1,p}}\leq c(\|D^{j}_{u}\xi\|_{L^{p}}+\|\xi\|_{L^{p}}) \end{equation} (for this
estimate, see Lemma B.4.6 in \cite{MS}, for example).  Where $\epsilon_{n}$, $\xi_{n}$ are as in the hypothesis, we have
\begin{align} \label{est}
\|\xi_{n}-\xi_{m}\|_{W^{1,p}}&\leq
c\big(\|D^{j}_{u}\xi_{n}-D^{j}_{u}\xi_{m}\|_{L^{p}}+\|\xi_{n}-\xi_{m}\|_{L^{p}}\big)
\notag
\\&=c\Big(\|(D^{j}_{u}-D^{j_{\epsilon_{n}}}_{u})\xi_{n}-(D^{j}_{u}-D^{j_{\epsilon_{m}}}_{u})\xi_{m}+D^{j_{\epsilon_{n}}}_{u}\xi_{n}-D^{j_{\epsilon_{m}}}_{u}\xi_{m}\|_{L^{p}}
\notag \\ & \qquad +\|\xi_{n}-\xi_{m}\|_{L^{p}} \Big) \notag \\
&\leq
c\big(C_{3}(\epsilon_{n}\|\xi_{n}\|_{W^{1,p}}+\epsilon_{m}\|\xi_{m}\|_{W^{1,p}})+\|D^{j_{\epsilon_{n}}}_{u}\xi_{n}\|_{L^{p}}+\|D^{j_{\epsilon_{m}}}_{u}\xi_{m}\|_{L^{p}}
 \notag \\ & \qquad  +\|\xi_{n}-\xi_{m}\|_{L^{p}}\big)
\end{align}
Now since $\{\xi_{n}\}$ is a bounded sequence in $W^{1,p}$, by Rellich compactness it has a subsequence which is Cauchy in $L^{p}$, and this fact along with the hypothesis of the lemma imply that, after
passing to a subsequence, the right hand side tends to zero as $m,n\to \infty$.  $\{\xi_{n}\}$ is therefore in fact Cauchy in $W^{1,p}$; say $\xi_{n}\to \xi$.  Then \[
D\xi=(D-D^{\epsilon_{n}})\xi + D^{\epsilon_{n}}(\xi-\xi_{n})+D^{\epsilon_{n}}\xi_{n}\to 0 \]
by Corollary \ref{lins} and the facts that $\xi_{n}\to \xi$ and $D^{\epsilon_{n}}\xi_{n}\to 0$.  But $D$ is injective, so $\xi=0$.  So the $\xi_{n}$ have a subsequence converging to zero.  If the
entire sequence did not converge to zero, we could take a subsequence bounded away from zero and apply the argument to that subsequence, obtaining a contradiction which proves the lemma. \end{proof}

\begin{cor}\label{inverses} \begin{itemize} \item[\rm(i)]  There is $\epsilon_{1}>0$ such that $D^{\epsilon}$ is bijective for all $\epsilon<\epsilon_1$.
\item[\rm(ii)] Denoting $Q^{\epsilon}=(D^{\epsilon})^{-1}$, for any sequence
$\epsilon_n \to 0$ we have $\|Q^{\epsilon_n}-Q\|\to 0$.
\end{itemize} \end{cor} \begin{proof}
If (i) were false we could find $\epsilon_n\to 0$ and $\xi_n$ with $\|\xi_n\|_{W^{1,p}}=1$ and $D^{\epsilon_n}\xi_n=0$.  This is prohibited by Lemma \ref{xin}.

For (ii), were this not the case for some sequence $\{\epsilon_n\}$, we could find
$\eta_n$ with $L^{p}$ norm 1 such that $Q^{\epsilon_n}\eta_n-Q\eta_n\nrightarrow 0$.  But then \begin{align}
\|D^{\epsilon_n}(Q^{\epsilon_n}\eta_n-Q\eta_n)\|_{L^{p}}&=\|D^{\epsilon_n}Q^{\epsilon_n}\eta_n+(D-D^{\epsilon_n})Q\eta_n-DQ\eta_n\|_{W^{1,p}} \notag \\
&=\|\eta_n+(D-D^{\epsilon_n})Q^{\eta_n}-\eta_n\|_{W^{1,p}}\leq C_{3}\|Q\|\epsilon_n \to 0 \notag \end{align}
violating Lemma \ref{xin} (with $\xi_n=Q^{\epsilon_n}\eta_n-Q\eta_n$) once again. \end{proof}

Corollary \ref{inverses} (ii) in particular implies that there is
$\epsilon_2 <\epsilon_1$ such that if $\epsilon < \epsilon_2$ then
$\|Q^{\epsilon}\|\leq \|Q\|+1$ (for otherwise we could find
$\epsilon_n\to 0$ with $\|Q^{\epsilon_n}-Q\|\geq 1$). Note that in
general, where $u\co (\Sigma,i)\to X$ denotes the (fixed)
embedding of $C$, we have
$\bar{\partial}_{j_{\epsilon}}u=\bar{\partial}_{j_{\epsilon}}u-\bar{\partial}_j
u=\frac{1}{2}(j_{\epsilon}-j)\circ du\circ i$, so since
$\|j_{\epsilon}-j\|_{C^{0}}\leq C_{2}\epsilon$ and
$j_{\epsilon}=j$ outside $D_{2r}^z \times D_{2\epsilon}^w$ (a
region whose intersection with $C$ has area proportional to
$\epsilon^2$), we have, for some constant $C_{4}$ related to $C_2$
and $\|du\|_{L^{\infty}}$, a bound
\begin{equation} \|\bar{\partial}_{j_{\epsilon}}u\|_{L^{p}}\leq C_{4}\epsilon^{1+2/p} \end{equation}
for $p>2$.  Fix such a $p$.  This puts us into position to prove:
\begin{lemma} \label{ift} There are constants $C_5$ and $\epsilon_3 >0$ such that for $\epsilon < \epsilon_3$ there exists
$\eta_{\epsilon}\in L^{p}(u^{*}N_C \otimes
T^{0,1}C)\oplus\bigoplus_{q\in\Omega}T_q X$ such that
$\bar{\partial}_{j_{\epsilon}}(\exp_u
(Q^{\epsilon}\eta_{\epsilon}))=0$ and
$\|Q^{\epsilon}\eta_{\epsilon}\|_{W^{1,p}}\leq
C_{5}\epsilon^{1+2/p}$.
\end{lemma}
\begin{proof}  This is a direct application of Theorem 3.3.4 of \cite{MS} (whose proof adapts without change to the case where the domain and range
consist of sections of $u^{*}N_C$ rather than $u^{*}TX$).  In
McDuff and Salamon's notation we take
$c_{0}=\max\{\|Q\|+1,\|du\|_{L^{p}},vol(\Sigma)\}$ and $\xi=0$.
The theorem gives $\delta$ and $c$ independent of $\epsilon$ such
that if $\|Q^{\epsilon}\|\leq c_{0}$ (as we have arranged to be
the case for $\epsilon < \epsilon_2$) and
$\|\bar{\partial}_{j_{\epsilon}}u\|_{L^{p}}\leq \delta$ then there
is $\eta_{\epsilon}$ with $\bar{\partial}_{j_{\epsilon}}(\exp_u
(Q^{\ep}\eta_{\epsilon}))=0$ and
$\|Q^{\epsilon}\eta_{\epsilon}\|\leq\|\bar{\partial}_{j_{\epsilon}}u\|_{L^{p}}$,
so we simply take $\epsilon_3 <\epsilon_2$ so small that
$C_{4}\epsilon_{3}^{1+2/p}\leq \delta$ and then $C_5 = cC_4$
\end{proof}

For $\epsilon < \epsilon_3$, let
$\xi_{\epsilon}=Q^{\epsilon}\eta_{\epsilon}$ and
$u_{\epsilon}=\exp_u\xi_{\epsilon}$.  We need to consider how the
branch points of the curve $C_{\ep}=u_{\ep}(\Sigma)$ relate to
those of $C$.  Our intent is to carry out this construction
sequentially for every branch point of $C$: at each step in the
procedure, then, we replace $j$ by an almost complex structure
which is integrable in some neighborhood of the branch point under
consideration, which has the effect of moving the curve somewhat;
we may assume inductively that at each of the previous steps our
procedure has resulted in the branch points being considered
getting replaced by branch points $p'$ contained in some
neighborhood $U'$ on which the new almost complex structure is
integrable.  For the present step, we need to ensure that two
things hold when $\ep$ is sufficiently small:
\begin{itemize} \item[(i)] That the branch points $q$ of $C_{\ep}$ that
are not close to $p$ are close enough to other branch points $p'$
of $C$ that if the neighborhood $U'$ as above (on which $j$ and so
also $j_{\ep}$ is integrable) has already been constructed around
$p'$, then $q\in U'$; and \item[(ii)] That the branch points of
$C_{\ep}$ which are close to $p$ fall into the neighborhood
$D_{r}^{z}\times D_{\ep}^{w}$ on which $j_{\ep}$ is integrable.
\end{itemize}

The first statement is somewhat easier, since every $j_{\epsilon}$
agrees with $j$ outside $D_{2r}^z \times D_{2\epsilon_0}^{w}$, and
so where $V$ is a small neighborhood of $D_{2r}^z \times
D_{2\epsilon_0}^{w}$  it follows from elliptic bootstrapping that
on $\Sigma\setminus u^{-1}(V)$ the $W^{1,p}$ bound on
$\xi_{\epsilon}$ implies $C^{k}$ bounds for all $k$.  Now all
branch points $p'$ of $C$ other than $p$ lie in $V$, so for any
such $p'$, since $f\circ u_{\ep}$ is holomorphic and tends to
$f\circ u$ in any $C^k$ norm near $p'$, for any neighborhood $U'$
of $u(p')$, if $\ep$ is small enough $U'$ will contain some number
$k$ of branch points $q_1,\ldots,q_k$ of $C_{\ep}$ such that,
where $n_q$ denotes the ramification index of a point $q$ on the
curve (equivalently, the order of tangency at $q$ between the
curve and the fiber), we have \[ \sum_m (n_{q_{m}}-1)=n_{p'}-1.\]
Conversely, at any $x\in\Sigma\setminus u^{-1}(V)$, the derivative
of $f\circ u_{\ep}$ at $x$ will be approximated to order
$\ep^{1+2/p}$ by that of $f\circ u$ at $x$.  In particular, if
$u_{\ep}(x)$ is a branch point, ie if $(f\circ
u_{\ep})_{*}$ is zero at $x$, then $(f\circ
u)_{*}(x)=O(\ep^{1+2/p})$, which if $\ep$ is small enough will
force $u(x)$ (and so also the new branch point $u_{\ep}(x)$, which
is a distance $O(\ep^{1+2/p})$ from $u(x)$) to be contained in any
previously-specified neighborhood of the branch locus of $C$. This
proves assertion (i) above.

Since the sum of the numbers $n_q-1$ where $q$ is a branch point
of $C_{\ep}$ is the same as the corresponding number for $C$ by
the Hurwitz formula applied to the holomorphic maps $f\circ
u_{\ep}$ and $f\circ u$, the sum of these numbers for just the
branch points of $C_{\ep}$ contained in $D_{2r}^z \times
D_{2\epsilon_0}^{w}$ must then $n_p-1$, $n_p$ being the
ramification index of $p$ as a branch point of $C$ (for by what
we've shown above, the sum of the $n_q-1$ for $q$ lying outside
this set also has not been changed by replacing $C$ with
$C_{\ep}$).

As such, $p$ is replaced either by a single branch point of
$C_{\ep}$ with ramification index $n_p$ or by some collection of
branch points (all in $D_{2r}^z \times D_{2\epsilon_0}^{w}$) each
with ramification index strictly less than $n_p$.  In the former
case, in the usual coordinates $(z,w)$ around $p$, since both $j$
and $j_{\ep}$ preserve all of the fibers $\{z=const\}$, as in
Section 2 of \cite{M} we may write $C$ as $\{z=w^{n_p}+O(n_p+1)\}$
and $C_{\ep}$ as $\{z=z_0 + k(w-w_0)^{n_p}+O(n_p+1)\}$ for some
$k$, where $(z_0,w_0)$ is the position of the new branch point.
But from Lemma \ref{ift} and the Sobolev Embedding theorem we have
an estimate $\|\xi_{\ep}\|_{C^{1-2/p}}\leq K\ep^{1+2/p}$, which
leads $z_0$, $k-1$, and $w_0$ to all be bounded by a constant
times $\ep^{1+2/p}$.  So if $\ep$ is small enough, the new node
$(z_0,w_0)$ will fall into the region $D^{z}_{r}\times
D_{\ep}^{w}$ on which $j_{\ep}$ is integrable, thanks to the fact
that $\ep^{1+2/p}\ll \ep$.

If instead $p$ is replaced by distinct branch points with lower
ramification indices, they in principle may not be so close, but
then we can apply our construction near each of these new branch
points.  Because at each step we either succeed or lower the
index, the process will eventually terminate (at the latest, when
the index has been lowered to two).

We should note that at each stage of the process the moduli space
only changes in the way that we have been anticipating.  Namely,
with the notation as above, we have:

\begin{lemma} \label{K}  Write $\mathcal{M}_{X}^{j,\Omega}(\alpha)=\{[ u ],[ v_1 ],\ldots ,[ v_r ]\}$.  Then for $\epsilon$ sufficiently small,
\[ \mathcal{M}^{j_{\epsilon},\Omega}_{X}(\alpha)=\{[ u_{\epsilon} ],[ v_1 ],\ldots ,[ v_r ]\}. \] \end{lemma}
\begin{proof}  That $\{[u_{\epsilon}],[v_1],\ldots ,[v_r]\}\subset \mathcal{M}^{j_{\epsilon},\Omega}_{X}(\alpha)$ is clear, since $u_{\epsilon}$ is
$j_{\epsilon}$--holo\-mor\-phic and passes through $\Omega$ by
construction (for it agrees with $u$ on the $u$--preimages of all
the points of $\Omega$), and since the $Im(v_{k})$ are all
contained in the set on which $j_{\epsilon}$ agrees with $j$.

To show the reverse inclusion, assume to the contrary that there
exists a sequence $\epsilon_n \to 0$ and $v_n\co \Sigma_n \to X$
with $[v_n]\in
\mathcal{M}^{j_{\epsilon_n},\Omega}_{X}(\alpha)\setminus
\{[u_{\epsilon_n}],[v_1],\ldots ,[v_r]\}$.  Now the almost complex
structures $j_{\epsilon_n}$ converge in the $C^{0}$ norm to $j$,
so by Gromov compactness (generalized to the case of $C^{0}$
convergence of the almost complex structures by Theorem 1 of
\cite{IS}), after passing to a subsequence there would be $[v]\in
\mathcal{M}_{X}^{j,\Omega}(\alpha)$ with $[v_{\epsilon_n}]\to [v]$
in any $W^{1,p}$ norm.  Now if $[v]$ were one of the $[v_k]$ this
would of course be impossible, since the $[v_{\epsilon_n}]$ would
then all eventually miss $D_{3r}^z \times D_{3\epsilon_0}^w$, so
the $Im(v_{\epsilon_n})$ would be contained in the region where
$j_{\epsilon_n}=j$, implying that the $v_{\epsilon_n}$ are
$j$--holomorphic curves passing through $\Omega$, which we assumed
they were not.

So suppose $[v_{\epsilon_n}]\to [u]$ in $C^{0}$.  Now
$u_{\epsilon_n}=\exp_u \xi_n$ with $\|\xi_n\|_{W^{1,p}}\leq C_5
\epsilon_{n}^{1+2/p}$, so
$\|u_{\epsilon_n}-v_{\epsilon_n}\|_{W^{1,p}} \to 0$ as $n\to
\infty$ for an appropriate parametrization of the
$v_{\epsilon_n}$. But, using the uniform boundedness of the right
inverses $Q^{\ep}$ of the linearizations $D^{j_{\epsilon}}_u$ at
$u$, Proposition 3.3.5 of \cite{MS} gives some $\delta$ such that
$\|u_{\epsilon_n}-v_{\epsilon_n}\|_{C^0}\geq \delta$ for all $n$,
a contradiction which proves the lemma.  \end{proof}

Lemma \ref{K} and the facts noted before it now let us prove the
following:

\begin{thm}  There is a constant $C_8$ such that for $\epsilon$ sufficiently small there exists an almost complex structure $\tilde{j}_{\epsilon}$ with $\|\tilde{j}_{\epsilon}-j\|_{C^0}\leq C_8 \epsilon$
having the property that, where
$\mathcal{M}^{\tilde{j}_{\epsilon},\Omega}_{X}(\alpha)=\{[
u_{1}^{\epsilon}], \ldots , [u_{r}^{\epsilon}]\}$,
$\tilde{j}_{\epsilon}$ is integrable on a neighborhood of each
point of $crit(f|_{Im(u_{i}^{\epsilon})})$.  Moreover
$\tilde{j}_{\epsilon}\in \mathcal{S}^0$, and
$\mathbb{J}_{\tilde{j}_{\epsilon}}$ is a regular almost complex
structure on $X_{r}(f)$. \end{thm}
\begin{proof}  Our construction shows how to modify $j$ into $j_{\epsilon}$ having the desired property in a small neighborhood of one branch point of one of the curves, say $C$, of $\mathcal{M}_{X}^{j,\Omega}(\alpha)$
without perturbing the other curves in
$\mathcal{M}_{X}^{j,\Omega}(\alpha)$, and, as noted above, the
construction can then be repeated at the other (slightly
perturbed) branch points of $C$, moving $C$ to a curve $C'$ near
\emph{all} of the branch points of which our new almost complex
structure has the desired property. Because the almost complex
structure remains unchanged near the other curves, we can apply
the same procedure sequentially to all of the curves of
$\mathcal{M}_{X}^{j,\Omega}(\alpha)$; this entails only finitely
many steps, at the end of which we obtain $\tilde{j}_{\epsilon}$,
which is regular by construction.

If $\mathbb{J}_{\tilde{j}_{\epsilon}}$ is not already regular,
Proposition \ref{JjOK} shows that it will become so after generic
perturbations of $\tilde{j}_{\epsilon}$ supported away from the
critical loci of the $f|_{Im(u_{i}^{\epsilon})}$ and the points of
$\Omega$. Provided they are small enough, such perturbations will
not change the other properties asserted in the theorem.
\end{proof}

\begin{cor} \label{good AC 2} In computing the invariant $Gr(\alpha)$, we can use an almost complex structure $j_1$ from the set $\mathcal{S}^1$ of Proposition \ref{JjOK}, and in computing the
invariant $\mathcal{DS}_{(X,f)}(\alpha)$, we can use the complex
structure $\mathbb{J}_{j_1}$.  \end{cor}

\section{Comparing the spectral flows}

We now fix an almost complex structure $j_1$ as in Corollary
\ref{good AC 2}, which we assume to have been constructed by the
procedure in the preceding section.  $C$ will now denote a fixed
member of $\mathcal{M}^{j_1,\Omega}_{X}(\alpha)$ with $u\co
(\Sigma,i)\to (X,j_1)$ a fixed embedding of $C$.  The assumption
on $\alpha$ at the start of the preceding section ensures that $C$
will not have any components which are multiply covered
square-zero tori; for more general $\alpha$ we now instead simply
assume that this is true for $C$.  We will show in this section
that the contribution of $C$ to $Gr(\alpha)$ is the same as that
of the associated section $s_C$ to $\mathcal{DS}_{(X,f)}(\alpha)$.

\begin{lemma} \label{intexists}  There is a neighborhood $U$ of $C$ and an integrable almost complex structure $\tilde{j}$ on $U$ which makes both $f$ and $C$ holomorphic.\end{lemma}

\begin{proof} Let $Crit(f|_C)=\{p_1,\ldots ,p_n\}$.  By our construction of $j_1$, there are neighborhoods $V_1,\ldots ,V_n$
of the $p_k$ on which $j_1$ is given by \[ T^{0,1}_{j_1}=\langle
\partial_{\bar{z}}+b(z,0)\partial_w, \partial_{\bar{w}} \rangle. \]
 Since all of the branch points of $C$ are contained within
 $\cup_k V_k$, we may cover $C\setminus \cup_k V_k$ by open sets
 $U_{\alpha}$ such that for each fiber $f^{-1}(t)$,
 $U_{\alpha}\cap f^{-1}(t)$ only contains at most one point of
 $C$.  In each $U_{\alpha}$, then, $C\cap U_{\alpha}$ is given as
 a graph
 \[ \{w_{\alpha}=\lambda_{\alpha}(z)\},\] where $w_{\alpha}$ is a
 $j_1$--holomorphic coordinate on the fibers; in such coordinates
 $C\cap U_{\alpha}$ will be holomorphic with respect to an almost
 complex structure given by $T^{0,1}=\langle \partial_{\bar{z}}+b(z,w_{\alpha})\partial_{w_{\alpha}}, \partial_{\bar{w}_{\alpha}} \rangle$ exactly if
 $b(z,\lambda_{\alpha}(z))=\frac{\partial\lambda}{\partial\bar{z}}$.
 We therefore simply define $\tilde{j}_{\alpha}$ on $U_{\alpha}$
 by \[
 T^{0,1}_{\tilde{j}_{\alpha}}=\langle \partial_{\bar{z}}+\frac{\partial\lambda}{\partial\bar{z}}\partial_{w_{\alpha}}, \partial_{\bar{w}_{\alpha}}
 \rangle.  \]  Geometrically, the $j_1|_{V_k}$ and the
 $\tilde{j}_{\alpha}$ are all uniquely determined by the fact that
 they restrict to the fibers as $j_1$, make $C$ and $f$
 holomorphic, and have defining functions $b$ which do not vary
 vertically, so in particular they agree on the overlaps of their
 domains and so piece together to form a complex structure
 $\tilde{j}$ on the set $U=\bigcup_k V_k\cup \bigcup_{\alpha}
 U_{\alpha}$, which is integrable by Equation \ref{Nj} and so
enjoys the properties stated in the lemma. \end{proof}

\begin{lemma} \label{sub} Let $\mathcal{J}(U,f,C)$ denote the set of almost complex structures on $U$ making both $C$ and $f$ holomorphic which are integrable near each branch point of $C$.  Let $\mathcal{J}^{int}(U,f,C)$ be the subset of $\mathcal{J}(U,f,C)$ consisting of almost complex structures integrable near all of $C$.
  Then the maps
\begin{align} \mathcal{F}\co H^{0,1}_i(T_{\mathbb{C}}\Sigma)\times W^{1,p}(u^{*}TX)\times \mathcal{J}(U,f,C) &\to L^p (u^{*}TX\otimes T^{0,1}\Sigma)   \notag \\
 (\beta,\xi,j)&\mapsto D_{u}^{j}\xi+\frac{1}{2}j\circ du\circ \beta, \notag \end{align} \begin{align}
\widehat{\mathcal{F}}\co W^{1,p}(s_{C}^{*}T^{vt}X_{r}(f)) \times \mathcal{J}(U,f,C) &\to L^{p}(s_{C}^{*}T^{vt}X_{r}(f)\otimes T^{0,1}S^2) \notag \\
(\zeta,j) &\mapsto D^{\mathbb{J}_j}_{s_C}\zeta \notag,
\end{align}
 \begin{align}  \mathcal{F}'\co H^{0,1}_i(T_{\mathbb{C}}\Sigma)\times W^{1,p}(u^{*}TX)\times \mathcal{J}^{int}(U,f,C) &\to L^p (u^{*}TX\otimes T^{0,1}\Sigma)   \notag \\
 (\beta,\xi,j)&\mapsto D_{u}^{j}\xi+\frac{1}{2}j\circ du\circ \beta, \notag \end{align} and
\begin{align}
\widehat{\mathcal{F}'}\co W^{1,p}(s_{C}^{*}T^{vt}X_{r}(f))\times \mathcal{J}^{int}(U,f,C) &\to L^{p}(s_{C}^{*}T^{vt}X_{r}(f)\otimes T^{0,1}S^2) \notag \\
(\zeta,j) &\mapsto D^{\mathbb{J}_j}_{s_C}\zeta \notag \end{align}
   are each submersive at all zeros whose section
component is not identically zero.
\end{lemma}
\begin{proof}
Suppose $\mathcal{F}(\beta,\xi,j)=0$. The linearization of
$\mathcal{F}$ at $(\beta,\xi ,j)$ is given by \begin{align}
\label{lin1} \mathcal{F}_{*}(\gamma,\mu, y)&=D^{j}_{u}\mu +
\left(\left.\frac{d}{dt}\right|_{t=0}
D_{u}^{exp_{j}(ty)}\right)\xi+\frac{1}{2}j\circ du\circ\gamma
\nonumber \\ &=D^{j}_{u}\mu + \frac{1}{2}(\nabla_{\xi}y)\circ
du\circ i+\frac{1}{2}j\circ du\circ\gamma, \end{align} where
$\nabla$ is the Levi--Civita connection of the metric associated
to $j$.   We assume $\xi$ is not identically zero, so that by
Aronzajn's theorem it does not vanish identically on any open
subset.
 If $\eta$ were a nonzero element of $\cok \mathcal{F}_{*}$, as
in the usual argument find $x_0\in\Sigma$ with  $u(x_0) \notin
Crit(f|_C)$ and $\eta(x_0)$ and $\xi(x_0)$ both nonzero. Near
$u(x_0)$, using the Levi--Civita connection of the metric
associated to $j$, $TX$ splits as $T^{vt}X\oplus TC$, and with
respect to this splitting $y$ (in order to be tangent to
$\mathcal{J}(U,f,C)$) is permitted to have any block decomposition
of form
\begin{equation}\label{y} y=\left( \begin{array}{cc} a & b \\ 0 & 0 \end{array} \right) \end{equation} where all entries are $j$--antilinear and, in order that $C$ remain holomorphic,
$b|_{C}=0$, so $\nabla_{\xi} y$ can have any block decomposition
of form $\left( \begin{array}{cc} a' & b'
\\ 0 & 0
\end{array} \right)$ where all entries are $j$--antilinear.
We have $0\neq \eta(x_{0})\in (u^{*}TX\otimes T^{0,1}\Sigma
)_{x_{0}}$, and $u(x_{0})\notin crit(f|_C)$, so
$(\eta(x_0))^{vt}\neq 0$. Hence similarly to the proof of Lemma
\ref{good AC} we can take $b'(x_0)$ and $c'(x_0)$ such that \[
\left(
\begin{array}{cc} 0 & b'(x_0) \\ 0 & 0 \end{array}
\right)du\circ i(v)=(\eta(x_0)(v))^{vt} \quad \left(
\begin{array}{cc} 0 & b'(x_0) \\ 0 & 0 \end{array}
\right)du\circ i(\bar{v})=(\eta(x_0)(\bar{v}))^{vt} \]  where $v$
generates $T^{1,0}_{x_0}\Sigma$.  We then take $y$ supported in a
small neighborhood of $u(x_0)$ so that $a=0$ in the decomposition
(\ref{y}) and so that
 \[\left(\nabla_{\xi}y\right)(x_0)=\left( \begin{array}{cc} 0 & b'(x_0) \\ 0 & 0 \end{array} \right) \]
By taking the small neighborhood appropriately, unless the
vertical projection $\eta^{vt}(x_0)$ of $\eta (x_0)$ is zero we
can thus arrange that
\[ \int \langle \eta,\mathcal{F}_{*}(0,y)\rangle \neq 0, \] in
contradiction with the assumption that $\eta$ belonged to the
cokernel of $\mathcal{F}_{*}$.  This shows that any $\eta\in \cok
\mathcal{F}_{*}$ must have $\eta^{vt}$ identically zero.  Then
arguing just as in the proof of Lemma \ref{good AC}, we consider
the projection $\eta^{C}$ of $\eta$ onto $TC$; once again $\eta^C$
would give an element of the cokernel of the linearization at
$(i,id)$ of the map $(i',v)\mapsto\dbar_{i',i}v$ acting on pairs
consisting of complex structures $i'$ on $\Sigma$ and maps $v\co
\Sigma\to\Sigma$, and the vanishing of this cokernel is just the
statement that the space of complex structures on $\Sigma$ is
unobstructed at $i$.  $\eta^C$ is therefore also zero, so since
$TC$ and $T^{vt}X$ span $TX$ at all but finitely many points of
$C$, we conclude that $\eta$ vanishes identically, proving the
Lemma for $\mathcal{F}$.

The proof of the transversality of $\widehat{\mathcal{F}}$
proceeds in essentially the same way; if $\eta\in \cok
(\widehat{\mathcal{F}}_{*})_{(\zeta,j)}$ with
$\widehat{\mathcal{F}}(\zeta,j)=0$ is nonzero at some $t$ (which
we can assume to be a regular value for $f|_C$), then as in the
proof of Lemma \ref{JjOK}, for at least one point $p_{0}$ among
the $r$ points of $X$ appearing in the divisor $s_{C}(t)$, $\eta$
descends to a nonzero $T_{p_0}^{vt}X$--valued cotangent vector at
$p_{0}$, and we can use a perturbation $y$ supported near $p_{0}$
similar to that above to obtain the desired contradiction.

As for $\mathcal{F}'$ and $\widehat{\mathcal{F}'}$, for which the
almost complex structure is required to be integrable near $C$,
the allowed perturbations $y$ include anything in the block form
\[ y=\left(
\begin{array}{cc} 0 & b \\ 0 & 0 \end{array} \right) \] where $b$
varies holomorphically in the vertical variable $w$ (as can be
seen from Equation \ref{Nj}).  So (aside from $j$--antilinearity)
we only require that for any vertical vector $\zeta$ we have
$\nabla_{j\zeta}b=j\nabla_{\zeta}b$.  For a particular tangent
vector $\xi$ at $u(x_0)$, then, we still have the freedom to make
$\nabla_{\xi}b$ any antilinear map that we choose, so we can just
duplicate the proof of the submersivity of $\mathcal{F}$ and
$\widehat{\mathcal{F}}$ to see that $\mathcal{F}'$ and
$\widehat{\mathcal{F}'}$ are also submersive at all zeros where
$\xi$ is not identically zero.
\end{proof}

\begin{cor} \label{goodintexists} There is a neighborhood $U$ of $C$ and an integrable almost complex structure $\tilde{j}$ on $U$ such that $\tilde{j}$ makes both $f$ and $C$
holomorphic, and such that the linearization
$\mathcal{D}^{\tilde{j}}_{u}$ of the operator
$(i,u)\mapsto\bar{\partial}_{i,\tilde{j}}u$ at the embedding of
$C$ is surjective, as is the linearization of
$\bar{\partial}_{\mathbb{J}_{\tilde{j}}}$ at $s_C$
\end{cor}
\begin{proof} We have just shown that the map $\mathcal{F}'\co H^{0,1}_i(T_{\mathbb{C}}\Sigma)\times(W^{1,p}(u^{*}TX)\setminus\{0\})\times \mathcal{J}^{int}(U,f,C) \to L^p (u^{*}TX\otimes
T^{0,1}\Sigma)$ which sends $(\beta,\xi,j)$ to
$\mathcal{D}^{j}_{u} (\beta,\xi)=D^{j}_{u}\xi+\frac{1}{2}j\circ
du\circ\beta$ is submersive at all zeros, so that the subset
$\{(\beta,\xi,j): \mathcal{D}^{j}_{u}(\beta,\xi)=0,\xi\not\equiv
0\}$ is a smooth manifold. As usual, applying the Sard--Smale
theorem to the projection onto the second factor we obtain that
for generic $j\in \mathcal{J}^{int}(U,f,C)$,
\[\ker\left( (\beta,\xi)\mapsto D^{j}_{u}\xi+\frac{1}{2}j\circ
du\circ\beta\right)\setminus
\{0\}=\ker\mathcal{D}^{j}_{u}\setminus\{0\}\] is a smooth manifold
of the expected dimension.  The correctness of the expected
dimension for generic $j\in \mathcal{J}^{int}(U,f,C)$ of course
translates directly to the surjectivity of the linearization
$\mathcal{D}^{j}_{u}$ for such $j$. Likewise, the submersivity of
$\widehat{\mathcal{F}'}$ shows that the linearization of
$\bar{\partial}_{\mathbb{J}_{\tilde{j}}}$ at $s_C$ is surjective
for generic $j\in\mathcal{J}^{int}(U,f,C)$. So since Lemma
\ref{intexists} shows that $\mathcal{J}^{int}(U,f,C)$ is nonempty,
the corollary follows.
\end{proof}

$\tilde{j}$ shall now denote an almost complex structure of the
type obtained by Corollary \ref{goodintexists}.

\begin{lemma} \label{paths}  There are paths $j_t$ of almost complex structures on $U$ connecting $j_0:=\tilde{j}$ to $j_1$ for which every $j_t$ makes both $f$ and $C$ holomorphic.  Moreover, for a dense set of
such paths: \begin{itemize} \item[\rm(i)] The path $j_t$ is
transverse to the set of almost complex structures $j$ for which
the linearization $D^{j}$ of the $\bar{\partial}_j$ operator at
$u$ (acting on normal sections) has excess kernel. \item[\rm(ii)] The
path $\mathbb{J}_{j_t}$ of complex structures on the subset
$\mathbb{U}$ of $X_{r}(f)$ corresponding to $U$ is transverse to
the set of almost complex structures $J$ for which the
linearization $D^J$ of the $\bar{\partial}_J$ operator at $s_C$
(acting on sections of $s_{C}^{*}T^{vt}X_{r}(f)$) has excess
kernel. \end{itemize}\end{lemma}

\begin{proof}  In local coordinates near $C$, the almost complex
structures $j_1$ and $\tilde{j}$ are given as \[
T^{0,1}_{j_1}=\langle
\partial_{\bar{z}}+b_1(z,w)\partial_w, \partial_{\bar{w}} \rangle
\] and \[T^{0,1}_{\tilde{j}}=\langle
\partial_{\bar{z}}+\tilde{b}(z,w)\partial_w, \partial_{\bar{w}}
\rangle.\]  Here we necessarily have $b_1|_C=\tilde{b}|_C$ since
both $j_1$ and $\tilde{j}$ make $C$ holomorphic, so to define a
path $j_t$ we can simply set \[ T^{0,1}_{j_t}=\langle
\partial_{\bar{z}}+((1-t)\tilde{b}(z,w)+tb_1 (z,w))\partial_w, \partial_{\bar{w}}
\rangle;\] on each chart (this obviously pieces together to give
an almost complex structure on all of $C$); since
$(1-t)\tilde{b}+tb_1 |_C=b_1 |_C =\tilde{b}|_C$, $C$ will be
$j_t$--holomorphic for each $t$.

As for statements (i) and (ii), Lemma \ref{sub} implies that the
map with domain \[ H^{0,1}_i(T_{\mathbb{C}}\Sigma)\times
(W^{1,p}(u^{*}N_C)\setminus\{0\})\times
(W^{1,p}(s_{C}^{*}T^{vt}X_{r}(f))\setminus\{0\}) \times
\mathcal{J}(U,f,C) \] defined by \[ (\beta,\xi,\zeta,j)\mapsto
(\mathcal{D}^{j}_{u}\xi, D^{\mathbb{J}_{j}}_{s_C}\zeta)
\] is transverse to zero, so that its zero set is a smooth
manifold and we obtain using the Sard--Smale theorem that, letting
$\mathbb{U}$ refer to the connected component containing $s_C$ in
the open subset of $X_{r}(f)$ consisting of unordered $r$--tuples
of points in $U\subset X$ that lie in the same fiber,
\begin{align*}
\mathcal{S}^1=\{j\in \mathcal{J}(U,f,C)\mid (j,\Omega),&
(\mathbb{J}_{j},\Omega)\\& \text{ are nondegenerate on $U$ and
$\mathbb{U}$ respectively}\}\end{align*} 
is open and dense; here
nondegeneracy of $(\mathbb{J}_{j},\Omega)$ means that the direct
sum $\mathbb{D}^j$ of $D^{\mathbb{J}_j}_{s_C}$ with the evaluation
map that tautologically corresponds to $(ev_{\Omega})_*$ is
bijective, while as in \cite{Taubes} nondegeneracy of $(j,\Omega)$
means that $D^{j}_{u}\oplus (ev_{\Omega})_*$ is bijective, which
is implied for generic $\Omega$ by the surjectivity of
$\mathcal{D}^{j}_{u}$.

Theorem 4.3.10 of \cite{DK} shows then that a dense set of paths
from $j_{0}$ to $j_1$ consists of paths which only cross the locus
for which either $D^{j}$ or $D^{\mathbb{J}_{j}}$ has excess kernel
transversely.  (Alternately we could of course prove a
parametrized version of Lemma \ref{sub} and apply the Sard--Smale
theorem to the projection to the space of paths in
$\mathcal{J}(U,f,C)$).\end{proof}
\begin{lemma} \label{kernels} For every $j\in \mathcal{J}(U,f,C)$ we have  \[ \ker (D^{j}_{u}\oplus (ev_{\Omega})_{*}) = 0 \iff \ker \mathbb{D}^j = 0.\] \end{lemma}
\begin{proof} Suppose that $\ker (D^{j}_{u}\oplus (ev_{\Omega})_{*}) \neq 0$ and let $0\neq \xi\in\ker (D^{j}_{u}\oplus (ev_{\Omega})_{*})$.
$\xi \in W^{1,p}\subset C^{0}$, so for $n$ sufficiently large $Im
(\exp_{u}(\xi / n))\subset U$.  Let $\eta_{n}$ be the sections of
$s_{C}^{*}T^{vt}X_{r}(f)$ such that $\exp_{s_{C}} \eta_n$
tautologically corresponds to $\exp_u (\xi/n)$.

By the construction of $\mathbb{J}_j$, for any point $t$ in the domain of $s_C$, $|\bar{\partial}_{\mathbb{J}_j}(\exp_{s_C}\eta_n)(t)|$ would be comparable
to the maximum of the  $|\bar{\partial}_j(\exp_u (\xi/n))|$ at the $r$ points corresponding to $s_C (t)$,
and similarly for $|\eta_n (t)|$ and the $|\xi /n|$ at the corresponding points, but for the fact that the end $q$ of a normal vector based at a point
$p_1 \in C$ will lie vertically over some other point $p_2 \in C$, which tends to increase distances as we pass to $X_{r}(f)$ since the (vertical) distance from
$p_2$ to $q$ will be larger than the length of the normal vector.  However, for any compact subset $K$ of $C\setminus crit(f|_C)$ normal vectors of small enough norm
based at some $p_1 \in K$ will correspond to vertical vectors based at some $p_2$ lying not too far outside of $K$ (and still outside of $crit (f|_C)$), and the norms
of the normal vector and the associated vertical vector will be comparable by some constant (depending on the set $K$).

Since as $n\to \infty$, $\exp_u (\xi/n)$ approaches the embedding
$u$ of $C$, we can then conclude the following: given $\epsilon$,
let $V_{\epsilon}\subset C$ be the $\epsilon$--neighborhood of
$crit(f|_C)$ in $C$.  Then there are $N$ and
$C_{1,\ep},C_{2,\ep},C_{3,\ep},C_{4,\ep}$ such that for $n\geq N$
we have:
\begin{equation} \label{etan} C_{1,\ep}\|\xi/n \|_{W^{1,p}(C\setminus V_{2\ep})} \leq \|\eta_n\|_{W^{1,p}(s_{C\setminus V_{\ep}})}\leq C_{2,\ep}\|\xi/n\|_{W^{1,p}(C\setminus V_{\ep/2})} \end{equation} and
\begin{align} \label{deletan} C_{3,\ep}\|\bar{\partial}_j \exp_u (\xi/n) \|_{L^{p}(C\setminus V_{2\ep})} &\leq \|\bar{\partial}_{\mathbb{J}_j} (\exp_{s_C} \eta_n)\|_{L^{p}(s_{C\setminus V_{\ep}})}\notag\\&\leq
C_{4,\ep}\|\exp_u (\xi/n)\|_{L^{p}(C\setminus V_{\ep/2})} \end{align}
Now since $D^{j} \xi =0$, there is a constant $C_5$ such that, for any $\ep, n$ we have \[
\|\bar{\partial}_j \exp_u (\xi/n) \|_{L^p (C\setminus V_{\ep})} \leq C_5 \|\xi /n \|_{W^{1,p}(C\setminus V_{\ep})}^{2} \]
Also, by Aronzajn's theorem, $\xi$ does not vanish on any open set, so writing $C_{6,\ep}=\frac{\|\xi\|_{W^{1,p}(C\setminus V_{\ep /2})}}{\|\xi\|_{W^{1,p}(C\setminus V_{2\ep})}}$, we have, independently
of $n$,
\[ \|\xi/n\|_{W^{1,p}(C\setminus V_{\ep /2})}\leq C_{6,\ep}\|\xi/n \|_{W^{1,p}(C\setminus V_{2\ep})} \]
We hence obtain, for all $n$ \begin{align}
\|\bar{\partial}_{\mathbb{J}_j}(\exp_{s_C} \eta_n)\|_{L^{p}s_(C\setminus V_{\ep})} &\leq C_{4,\ep}\|\exp_u (\xi/n)\|_{L^{p}(C\setminus V_{\ep/2})}\notag\\
&\leq C_{4,\ep}C_5  \|\xi /n \|_{W^{1,p}(C\setminus V_{\ep/2})}\notag\\
& \leq C_{4,\ep}C_5 C_{6,\ep}^{2} \|\xi/n \|_{W^{1,p}(C\setminus V_{2\ep})}^{2}\notag\\ &\leq \frac{C_{4,\ep}C_5 C_{6,\ep}^{2}}{C_{1,\ep}^{2}}  \|\eta_n\|_{W^{1,p}(s_{C\setminus V_{\ep}})}^{2}\notag \end{align}
So we have $W^{1,p}$ sections $\eta_n \to 0$ of
$s_{C}^{*}T^{vt}X_{r}(f)$ such that, for each $\epsilon$,
\begin{equation}\label{tozero}
\frac{\|\bar{\partial}_{\mathbb{J}_j}(\exp_{s_C}
\eta_n)\|_{L^{p}(s_{C\setminus
V_{\ep}})}}{\|\eta_n\|_{W^{1,p}(s_{C\setminus V_{\ep}})}}\to 0
\end{equation}
We now show how to obtain from (\ref{tozero}) an element of the kernel of the linearization $D^{\mathbb{J}_j}_{s_C}$.

Fix $\ep$ and consider the linearization $D_{\ep}$ of
$\bar{\partial}_{\mathbb{J}_j}$ at $s_{C\setminus V_{\ep}}$,
acting on $W^{1,p}$ sections of the bundle $E_{\ep}=s_{C\setminus
V_{\ep}}^{*}T^{vt}X_{r}(f)$.  Let $r_{n}\co E_{\ep}\to E_{\ep}$ be
the bundle endomorphism given by fiberwise multiplication by
$\frac{1}{\|\eta_n\|_{W^{1,p}(s_{C\setminus V_{\ep}})}}$.
Identifying a neighborhood of the zero section in $E_{\ep}$ with a
neighborhood of $s_{C\setminus V_{\ep}}$, we have that, fixing $k$
small enough that each $Im\left(\exp_{s_{C\setminus
V_{\ep}}}\left(\frac{k\eta_n}{\|\eta_n\|_{W^{1,p}(s_{C\setminus
V_{\ep}})}}\right)\right)$ is in this neighborhood (which is
possible since the $\eta_n /\|\eta_n \|$ are $C^{0}$--bounded),
\[
\bar{\partial}_{r_{n}^{*}\mathbb{J}_{j}}\left(\exp_{s_{C\setminus
V_{\ep}}}\left(\frac{k\eta_n}{\|\eta_n\|_{W^{1,p}(s_{C\setminus
V_{\ep}})}}\right)\right)
=\frac{k}{\|\eta_n\|_{W^{1,p}(s_{C\setminus
V_{\ep}})}}\bar{\partial}_{\mathbb{J}_j}\exp_{s_{C\setminus
V_{\ep}}}\eta_n \to 0,\] and each
$\frac{k\eta_n}{\|\eta_n\|_{W^{1,p}(s_{C\setminus V_{\ep}})}}$ has
norm $k$.  Write $\zeta_n =
\frac{k\eta_n}{\|\eta_n\|_{W^{1,p}(s_{C\setminus V_{\ep}})}}$.

Now since $r_{n}$ is multiplication by  $\frac{1}{\|\eta_n\|_{W^{1,p}(s_{C\setminus V_{\ep}})}}$, which tends to $\infty$ with $n$, we have that \[
\lim_{n\to \infty} D_{\ep}\zeta_n =  \lim_{n\to \infty}\bar{\partial}_{r_{n}^{*}\mathbb{J}_j}(\exp_{s_{C\setminus V_{\ep}}}\zeta_n) =0 \]
By Rellich compactness, after passing to a subsequence the
$\zeta_n$ $L^p$--converge to some $\zeta^{\ep}\in L^{p}$; since
the $\zeta_n$ have norm bounded away from zero, $\zeta^{\ep}\neq
0$.  Where $D_{\ep}^{*}$ is the formal adjoint of $D_{\ep}$, we
then have that, for each $\beta\in W^{1,q}(\Lambda^{0,1}M_P\otimes
s_{C\setminus V_{\ep}}^{*}T^{vt}X_{r}(f))$ ($1/p+1/q=1$),
\[ \langle \zeta^{\ep}, D^{*}_{\ep}\beta \rangle=\lim_{n\to \infty}\langle \zeta_n,D^{*}_{\ep}\beta\rangle = \lim_{n\to \infty}\langle D_{\ep}\zeta_n , \beta \rangle =0\]
So $\zeta^{\ep}$ is a weak solution to $D_{\ep}\zeta_{\ep}=0$; by elliptic regularity this implies that $\zeta^{\ep}$ is in fact in $W^{1,p}$ with $D_{\ep}\zeta^{\ep}=0$.

All of the $\zeta^{\epsilon}$ so constructed agree up to scale on
the overlaps of their domains (since they are limits of rescaled
versions of the $\eta_n$, and the $\eta_n$ do not vary with
$\epsilon$); also if we require that the tubular neighborhoods of
$s_{C\setminus V_{\ep}}$ used in the construction are all
contained in a common tubular neighborhood of $s_C$, the
$\exp_{s_C} \zeta^{\epsilon}$ will all be contained in this
neighborhood, so that the norms of the $\zeta^{\ep}$ will be
bounded, say by $M$, as $\ep \to 0$.  So we can rescale the
$\zeta^{\ep}$ to all agree on their domains with a common section
$\zeta\in W^{1,p}(s_{C}^{*}T^{vt} X_{r}(f))$ defined on the
complement from the finite set of critical values of $f|_C$ which
is nonzero (since all of the $\zeta^{\epsilon}$ are) and has
$D_{\ep}\zeta =0$ for every $\epsilon>0$.  Moreover the norm of
$\zeta$ on any compact subset of its domain is at most $M$, so by
removal of singularities $\zeta$ extends to all of $S^2$, and
$\zeta \in \ker D^{\mathbb{J}_j}_{s_C}$. Further, since $\xi \in
\ker (ev_\Omega)_{*}$, it readily follows from the construction
that $\zeta$ is in the kernel of the corresponding linearization
of the corresponding evaluation map on $X_{r}(f)$, so that $0\neq
\zeta\in \ker \mathbb{D}^j$, proving the forward implication in
the statement of Lemma \ref{kernels}.

The reverse implication can be proven in just the same way, by
taking an element $0\neq \eta \in \ker \mathbb{D}^j$ and
extracting a normal section $\xi$ from the curves tautologically
corresponding to the $\exp (\eta/n)$ which lies in the kernel of
the restriction of $(D^{j}_{u}\oplus (ev_{\Omega})_{*})$ to any
set missing $crit(f|_C)$  Once again, removal of singularities
then implies that $\xi$ extends to give a global nonzero element
of $\ker (D^{j}_{u}\oplus (ev_{\Omega})_{*})$.
\end{proof}

This directly yields the theorem promised at the beginning of the
section.
\begin{thm} \label{singleagree}  The contribution of $C$ to $Gr(\alpha)$ is the same as that of $s_{C}$ to $\mathcal{DS}_{(X,f)}(\alpha)$.  \end{thm}
\begin{proof}
Take a path $j_{t}$ as in Lemma \ref{paths}, so that $j_t$ is
transverse to the set of $j$ for which either $D^{j}_{u}\oplus
(ev_{\Omega})_{*}$ or $\mathbb{D}^{j}$ has nonzero kernel.  Since
$N_{\tilde{j}}=0$, we have $N_{\mathbb{J}_{\tilde{j}}}=0$, so by
the remarks at the start of Section \ref{gac2} the contribution of
$C$ to $Gr$ may be computed from the spectral flow of the path of
operators $D^{j_t}_{u}\oplus (ev_{\Omega})_{*}$, while that of
$s_C$ to $\mathcal{DS}$ may be computed from the spectral flow of
the path $\mathbb{D}^{j_t}$. By Lemma \ref{kernels}, for every $t$
the operator $D^{j_t}_{u}\oplus (ev_{\Omega})_{*}$ has a kernel if
and only if $\mathbb{D}^{j_t}$ does, so the number of eigenvalue
crossings for positive $t$, each of which is known to be
transverse, will be the same.  The two contributions are then both
equal to negative one to this common number of crossings.
\end{proof}

\section{Multiple covers of square-zero tori}

For curves with square-zero toroidal components, the difficulties
involved in comparing the contributions to $Gr$ and $\mathcal{DS}$
are more serious.  On the $Gr$ side, as Taubes showed in
\cite{Taubes}, if $C$ is a $j$--holomorphic square-zero torus, not
only $C$ but also each of its multiple covers contributes to $Gr$,
according to a prescription which depends on the spectral flows
not only of the linearization $D$ of the $\dbar$ operator on the
normal bundle $N_C$ but also of the three operators $D_{\iota}$
corresponding to $D$ which act on sections of the bundle obtained
by twisting $N_C$ by the real line bundles with Stiefel--Whitney
class $\iota$.  From the standpoint of the tautological
correspondence, it is encouraging that multiple covers of
square-zero tori contribute to $Gr$, since such covers do
tautologically correspond to $\mathbb{J}_j$--holomorphic sections
of $X_{r}(f)$ for appropriate $r$.  These sections are more
difficult to analyze, though, because they are contained in the
diagonal stratum $\Delta$ of $X_{r}(f)$, so the problems stemming
from the nondifferentiability of $\mathbb{J}_j$ cannot be evaded
by modifying $j$ to be integrable near the branch points.

Throughout this section, all almost complex structures $j$ defined on some region of $X$
that we consider will be assumed to make the restriction of $f$ to that region pseudoholomorphic.

As in Definition 4.1 of \cite{Taubes}, a $j$--holomorphic
square-zero torus $C$ will be called $m$--nondegenerate if, for
each holomorphic cover $\tilde{C}\to C$ of degree at most $m$, the
operator $\tilde{D}$ obtained by pulling back the linearization
$D$ (which acts on $\Gamma (u^* N_C)$ if $u$ is the map of $C$
into $X$) by the cover $\tilde{C}\to C$ has trivial kernel.  $j$
will be called $m$--nondegenerate for some fixed cohomology class
$\alpha\in H^2 (X,\mathbb{Z})$ with $\alpha^2 = \kappa\cdot\alpha
=0$ if every $j$--holomorphic curve $C$ with $[C]=PD(\alpha)$ is
$m$--nondegenerate. Lemma 5.4 of \cite{Taubes} shows that
$m$--nondegeneracy is an open and dense condition on $j$.

For any integer $m$, if $C$ is a $j$--holomorphic square-zero
torus Poincar\'e dual to the class $\alpha$, where  $j$ is
$m$--nondegenerate and is as in Lemma \ref{good AC}, we can define
the contribution $r'_{j}(C,m)$ of $m$--fold covers of $C$ to
$\mathcal{DS}_{(X,f)}(m\alpha)$ as follows.  Take a small tubular
neighborhood $U$ of $C$ which does not meet any of the other
$j$--holomorphic curves Poincar\'e dual to any  $k\alpha$ where
$k\leq m$ (this is possible since the nondegeneracy of $j$ ensures
that there are only finitely many such curves and since
$\alpha^2=0$) and which misses the critical points of the
fibration. Where $r$ is the intersection number with the fibers of
$f$, let $\mathbb{U}$ be the neighborhood of the section $s_{mC}$
of $X_{mr}(f)$ tautologically corresponding to $U$, so
$\mathbb{J}_j$ is H\"older continuous (say $C^{\gamma}$) on
$\mathbb{U}$ and $s_{mC}$ is the only $\mathbb{J}_j$--holomorphic
section in its homotopy class which meets $\mathbb{U}$.  Let $V$
be an open set with closure contained in $\mathbb{U}$ and
containing the image of $s_{mC}$; then it follows readily from
Gromov compactness that there is $\ep>0$ such that if $J$ is any
almost complex structure with
$\|J-\mathbb{J}_j\|_{C^{\gamma}}<\ep$ then any $J$--holomorphic
curve meeting $\mathbb{U}$ must in fact be contained in $V$.
$r'_{j}(C,m)$ is then defined as the usual signed count of all
$J$--holomorphic sections homotopic to $s_{mC}$ and contained in
$V$ where $J$ is a generic almost complex structure which is
smooth on $V$ and has $\|J-\mathbb{J}_j\|_{C^{\gamma}}<\ep$.  The
usual cobordism argument (using cobordisms which stay
H\"older-close to $\mathbb{J}_j$ so that sections in the
parametrized moduli spaces don't wander outside of $V$) shows that
this count is independent of the choice of $J$. Similarly, for any
$\beta\in H^{2}(X,\mathbb{Z})$, defining the contribution to
$\mathcal{DS}_{(X,f)}(\beta)$ of any disjoint union of
$j$--holomorphic curves with multiplicities with homology classes
adding to $PD(\beta)$ by smoothing $\mathbb{J}_j$ near the
associated section of $X_{r}(f)$, one notes that
$\mathcal{DS}_{(X,f)}(\beta)$ is indeed the sum of all the
contributions of all such unions, so the terminology is not
misleading.

Note that this definition of the contribution of $m$--fold covers
of $C$ to $\mathcal{DS}$ makes sense even if $C$ is itself a
multiple cover.  If $C$ is a $k$--fold cover of $C'$, then the
section $s_{lC}$ associated to an $l$--fold cover of $C$ is just
the same as the section $s_{klC'}$, and $r'_j(C,l)$ is defined by
perturbing the almost complex structure on the relative Hilbert
scheme near this section.  In particular, we have
$r'_j(C,l)=r'_j(C',kl)$.

\begin{lemma} \label{nowalls1} Let $j_{t}$ $(0\leq t\leq 1)$ be a path of almost complex structures which
make $f$ holomorphic such that every $j_{t}$ is
$m$--non-degenerate, and let $C_{t}$ be a path of embedded
square-zero tori in $X$ such that $\{(C_t,t)|0\leq t\leq 1\}$ is
one of the connected components of the parametrized moduli space
of $j_t$--holomorphic curves homologous to $C_0$. Then
$r'_{j_0}(C_0,m)=r'_{j_1}(C_1,m)$.
\end{lemma}

\begin{proof} Because all of the $j_t$ are $m$--non-degenerate,
there is an open neighborhood $U$ of $\cup_t C_t\times\{t\}\subset
X\times [0,1]$ such that no curve in homology class $k[C_t]$ for
any $k\leq m$ meets $U$ (for otherwise Gromov compactness would
give either a $j_t$--holomorphic curve in class $k[C_t]$ meeting
$C_t$ in an isolated point, which is impossible since $[C_t]^2=0$,
or a sequence of curves distinct from $C_t$ which converge to a
$k$--fold cover of $C_t$, which is prohibited by
$m$--non-degeneracy). Where $r$ is the intersection number of
$C_t$ with the fibers of $f$, let $\mathbb{U}$ be the neighborhood
of $\cup_t Im(s_{mC})\times\{t\}$ tautologically corresponding to
$U$ and $V$ some neighborhood of $\cup_t Im(s_{mC})\times \{t\}$
compactly contained in $\mathbb{U}$. Let $J_t$ be a family of
smooth almost complex structures on $X_{mr}(f)$ which are
sufficiently H\"older-close to $\mathbb{J}_{j_t}$ that each
$J_t$--holomorphic section meeting $\mathbb{U}$ is contained in
$V$, taken so that $J_0$ and $J_1$ are both regular and the path
$J_t$ is suitably generic. Now $\{(s,t)|\dbar_{J_t}s=0\}$ of
course gives an oriented cobordism between the moduli spaces of
$J_0$ and $J_1$--holomorphic sections in the relevant homotopy
class, and moreover, since none of the members of
$\{(s,t)|\dbar_{J_t}s=0\}$ even meet the open set
$\mathbb{U}\setminus \bar{V}$, this cobordism restricts to a
cobordism between the set of $J_0$--sections contained in $V$ and
the set of $J_1$--sections contained in $V$.  Since the
$r'_{j_k}(C_k,m)$ ($k=0,1$) are precisely the signed count of
these sections, it follows that $r'_{j_0}(C_0,m)=r'_{j_1}(C_1,m)$.
\end{proof}

A major reason that the analysis of multiply-covered
pseudoholomorphic curves is generally more difficult is that when
multiply-covered curves are allowed the argument that is generally
used to show the submersivity of the ``universal map'' $(u,j)
\mapsto \bar{\partial}_j u$ breaks down.  As a consequence, for
instance, as far as the author can tell it is not possible to
ensure that a square-zero torus $C$ will admit any almost complex
structures near it which both make it $m$--nondegenerate and are
integrable if $m>1$.  In the semi-positive context in which we
presently find ourselves, the standard way to navigate around this
difficulty, following \cite{RT1} and \cite{RT}, is to construct
our invariants from solutions to the inhomogeneous Cauchy--Riemann
equation \begin{equation} \label{inhomog} (\bar{\partial}_j
u)(p)=\nu(p, u(p)), \end{equation} where the domain of the map
$u\co \Sigma\to X$ is viewed as contained in a ``good cover'' of
the universal curve $\bar{\mathcal{U}}_{g,n}$ which is itself
embedded in some $\mathbb{P}^{N}$, and $\nu$ is a section of the
bundle $Hom (\pi_{1}^{*}T\mathbb{P}^{N},\pi_{2}^{*}TX)\to
\mathbb{P}^N \times X$ which is antilinear with respect to the
standard complex structure on $\mathbb{P}^N$ and the almost
complex structure $j$ on $X$ (see Definitions 2.1 and 2.2 of
\cite{RT} for details; note however in our case since we are
counting curves which may not be connected, we need to replace
$\bar{\mathcal{U}}_{g,n}$ with the universal space
$\bar{\mathcal{U}}^{(m)}_{\chi,n}$ of curves with at most $m$
components, $n$ marked points, and total Euler characteristic
$\chi$). Solutions to this equation are called
$(j,\nu)$--holomorphic curves. $\nu$ is called an inhomogeneous
term.

Imitating very closely the proof of Lemma \ref{sub}, one can see
that for any given $m\geq 1$ and for any fixed
$(j,0)$--holomorphic curve $C$ and for generic inhomogeneous terms
$\nu$ which
\begin{itemize} \item[(a)] vanish along the graphs of the embedding
$u$ of $C$ and of all of its covers up to degree $m$, \item[(b)]
take values in $T^{vt}X$ (rather than just $TX$), \item[(c)] are
``holomorphic in the X variable'' in the sense that
$\nabla_{(0,j\zeta)}\nu = j\nabla_{(0,\zeta)}\nu$ for $\zeta\in
TX$ (and  $(0,\zeta)\in T(\mathbb{P}^N\times X)$), and \item[(d)]
have the following ``coherence'' property: where $u\co \Sigma\to
X$ is embedding of $C$ and $\phi'\co \Sigma'\to \Sigma$ and
$\phi''\co \Sigma''\to \Sigma$ are any two holomorphic, possibly
disconnected, $m$--fold covers of $\Sigma$, for each $p \in
\Sigma$ and each $x\in X$ close to $u(p)$ the unordered
$m$--tuples $\{\nu(p',x):\phi'(p')=p\}$ and
$\{\nu(p'',x):\phi''(p'')=p\}$ are the same,
\end{itemize}
all of the covers of $C$ of degree $m$ will be nondegenerate as
$(j,\nu)$--holomorphic curves (ie, the linearization of
the equation (\ref{inhomog}) will be surjective at each of these
covers).  The point of condition (c) above is that it ensures that
these linearizations are all complex linear if $j$ is integrable
near $C$.  The point of condition (d) is that it ensures that
there is an inhomogeneous term $\mu$ on $X_{mr}(f)$ such that the
equation for a $(j,\nu)$--holomorphic curve in class $m[C]$ near
$C$ is the same as the equation for a
$(\mathbb{J}_j,\mu)$--holomorphic section of $X_{mr}(f)$ near
$s_{mC}$ which descends to a cycle in class $m[C]$.  $\nu$
satisfying this condition may easily be constructed: any choice of
$m$ perturbation terms
$\nu_1,\ldots,\nu_m\in\Gamma(\overline{Hom}(T\Sigma,u^{*}T^{vt}X))$
which vanish near the branch points of $C$ can be assembled into
perturbation terms near each of the holomorphic $m$--fold covers,
and we can use cutoff functions to put these together in order to
form a coherent inhomogeneous term $\nu\in
\Gamma(\overline{Hom}(\pi^{*}_{1} T\mathbb{P}^N,\pi^{*}_{2}TX))$.
Since the curves giving $m$--fold covers of $\Sigma$ in
$\bar{\mathcal{U}}^{(m)}_{\chi=0,n}$ are separated from each
other, the coherence condition does not make the proof of generic
nondegeneracy any more difficult. The reason that we can imitate
the proof of Lemma \ref{sub} using inhomogeneous terms but not
using almost complex structures is of course that we need the
freedom to vary the linearization of the equation on individual
small neighborhoods in the domain while leaving it unchanged
elsewhere, and for, say, a $k$--fold cover, varying the almost
complex structure on a small neighborhood in $X$ has the effect of
varying the linearization on $k$ different neighborhoods of the
domain all in the same way.

A pair $(j,\nu)$ such that $\nu$ satisfies conditions (b) through
(d) with respect to all $(j,\nu)$--holomorphic curves $C$ will be
called \emph{admissible}.  We will slightly enlarge the class of
data we study as follows: instead of only considering pairs
$(C,j)$ where $C$ is $j$--holomorphic, we consider triples
$(C,j,\nu)$ where $C$ is $j$--holomorphic, $\nu$ vanishes along
the graphs of the embedding of $C$ and all of its covers up to
degree $m$, and $(j,\nu)$ is admissible; such a triple will be
called $m$--nondegenerate if all of the covers of $C$ of degree
$m$ or lower are nondegenerate as $(j,\nu)$--holomorphic curves.
The admissible pair $(j,\nu)$ will itself be called
$m$--nondegenerate if $(C,j,\nu)$ is $m$--nondegenerate for each
$(j,\nu)$--holomorphic curve $C$. We can then define the
contribution $r'_{j,\nu}(C,m)$ to $\mathcal{DS}$ if $(C,j,\nu)$ is
$m$--nondegenerate: the nondegeneracy implies that there is a
neighborhood $U$ of $C$ which does not meet any other
$(j,\nu)$--holomorphic curves in class $k[C]$ for $k\leq m$.  We
have a tautologically-corresponding inhomogeneous term $\mu$ on
$X_{mr}(f)$, and we may perturb the almost complex structure
$\mathbb{J}_j$ to a smooth almost complex structure $J$ such that
$(J,\mu)$ is nondegenerate on a neighborhood $V$ of $s_{mC}$
contained in the set tautologically corresponding to $U$; we then
count $(J,\mu)$ holomorphic sections according to the prescription
in \cite{RT1}.  (Gromov compactness in the context of solutions to
the inhomogeneous Cauchy--Riemann equation is needed here; this
result appears as Proposition 3.1 of \cite{RT1}.)  The proof of
Lemma \ref{nowalls1} then goes through to show:

\begin{cor} \label{nowalls} Let $(j_{t},\nu_t)$ $(0\leq t\leq 1)$
be a path of $m$--nondegenerate admissible pairs, and let $C_{t}$
be a path of embedded square-zero tori in $X$ such that
$\{(C_t,t)|0\leq t\leq 1\}$ is one of the connected components of
the parametrized moduli space of $(j_t,\nu_t)$--holomorphic curves
homologous to $C_0$. Then
$r'_{j_0,\nu_0}(C_0,m)=r'_{j_1,\nu_1}(C_1,m)$.
\end{cor}

Now assume that $(C,j,\nu)$ is $m$--nondegenerate and that $j$ is
\emph{integrable} near $C$.  $\mathbb{J}_j$ is then smooth (and
even integrable) near $s_{mC}$; the argument in the proof of Lemma
\ref{kernels} shows that $(\mathbb{J}_j,\mu)$ will then also be
nondegenerate (and even if it weren't, it would become so after a
suitable perturbation of $\nu$ among inhomogeneous terms
satisfying conditions (a) through (d)), so in computing
$r'_{j,\nu}(C,m)$ we don't need to perturb $\mathbb{J}_j$ at all.
So since the linearization of the equation
$\bar{\partial}_{\mathbb{J}_j}s=\mu$ at $s_{mC}$ is complex-linear
and since $s_{mC}$ is the only solution to that equation in $V$,
we obtain (using Corollary \ref{nowalls}):

\begin{lemma} \label{intcase} If $(j,\nu)$ is an admissible pair and $C$ a $j$--holomorphic square-zero torus
such that $j$ is integrable near $C$, and if the
$m$--non-degenerate pair $(j',\nu')$ with $C$ $j'$--holomorphic is
sufficiently close to $j$, then $r'_{j',\nu'}(C,m)=1$ for every
$m$.
\end{lemma}

Our basic strategy in proving that multiple covers of square-zero
tori contribute identically to $\mathcal{DS}$ and $Gr$ will be,
using an almost complex structure $j$ as in Corollary \ref{good AC
2}, to investigate how the contributions $r'_{j_t,\nu_t}(C,m)$
vary as we move among admissible pairs such that $C$ is
$j_t$--holomorphic along a path from an $m$--nondegenerate pair
$(j_0,\nu_0)$ with $j_0$ integrable near $C$ to the pair $(j,0)$
where $j$ is the given nondegenerate almost complex structure.
This requires a digression into the chamber structure of almost
complex structures on $X$, which was investigated extensively by
Taubes in \cite{Taubes}.  For simplicity of exposition, we will
generally work in the homogeneous context $\nu=0$ below; since the
wall crossing results that follow only depend on the basic shape
of the differential equations involved and their linearizations,
the results below will remain valid when ``$j_t$'' is replaced by
$``(j_t,\nu_t).$''

Where $\mathcal{M}_{1,1}$ is the moduli space of smooth pointed
complex tori, consider the bundle $\mathcal{G}\to
\mathcal{M}_{1,1}$ whose fiber over the curve $C$ is the set of
1--jets at $C$ of almost complex structures on the trivial complex
line bundle over $C$.  Any such 1--jet gives rise to four
linearizations $D_{\iota}$ of the $\dbar$ operator on the bundles
$\underline{\mathbb{C}}\otimes L_{\iota}$ over $C$, where
$L_{\iota}$ is the real line bundle over $C$ with Stiefel--Whitney
class $\iota \in H^{1}(C,\mathbb{Z}/2)$.  Taubes shows that the
set $\mathcal{D}_{\iota}$ of points of $\mathcal{G}$ whose
corresponding linearization has a nontrivial kernel is a
subvariety of real codimension at least 1, and that the set of
elements of $\mathcal{D}_{\iota}$ either corresponding to a
linearization with two-or-greater-dimensional kernel or belonging
to some other $\mathcal{D}_{\iota'}$ has real codimension at least
2 in $\mathcal{G}$.  Identical results apply when we instead take
the fiber of $\mathcal{G}$ to consist of 1--jets of admissible
pairs $(j,\nu)$.

A path $\gamma=(u_t,C_t,j_t)_{t\in [0,1]}$ of $j_{t}$--holomorphic
immersions $u_{t}\co C_t \to X$ (each $C_t$ belonging to
$\mathcal{M}_{1,1}$; more commonly we will just denote such paths
by $(C_t, j_t)$, suppressing the map and identifying $C_t$ with
its image in $X$) then gives rise to a path $\tilde{\gamma}$ in
$\mathcal{G}$; we say $\gamma$ crosses a wall at $t=t_0$ if
$\tilde{\gamma}$ meets one of the codimension-one sets
$\mathcal{D}_{\iota}$ transversely at $t_0$. (Note that it's not
essential that the $u_t$ be embeddings, and in fact the case where
$u_t$ is a double cover will be of some relevance later on). The
path components of $\mathcal{G}\setminus
\cup_{\iota}\mathcal{D}_{\iota}$ are called \emph{chambers}. For
any $m$, Part 5 of Lemma 5.8, Lemma 5.9, and Lemma 5.10 of
\cite{Taubes} show (among other things) that for a generic path
$(C_t, j_t)$, the only $t_0$ for which  $j_{t_0}$ fails to be
$m$--nondegenerate near $C_{t_0}$ are those $t_0$ for which
$(C_{t_0}, j_{t_0})$ is on a wall. The proofs of the results
concerning connectivity and regularity of almost complex
structures which make $f$ holomorphic from sections 2 through 4
may easily be modified to show that the corresponding statement is
true for paths $j_t$ generic among paths of almost complex
structures which make $f$ holomorphic. On a similar note, if a
path $(C_t,j_t)$, where each $j_t$ is an almost complex structure
which makes $f$ holomorphic, remains in the same chamber except
for one point at which it touches a wall, the arguments in the
proofs of Lemmas \ref{good AC} and \ref{paths} show that the path
may be perturbed to a path which remains entirely within the
chamber and for which the almost complex structure continues to
make $f$ holomorphic.

In general, with the convention that $r'_{j}(C,0)=1$, we will
organize the contributions $r'_j (C,m)$ into a generating function
(to be viewed as a formal power series; we are not making any
convergence assertions here)
\[ P'_j (C,z)=\sum_{m\geq 0} r'_j (C,m)z^m.\]
Strictly speaking, this power series should be truncated after the
term corresponding to the largest $m$ for which $j$ is
$m$--non-degenerate and the fibration satisfies $\omega\cdot
(fiber)
>m\omega\cdot \alpha$.  However, by working with suitably generic $j$
and suitably high-degree Lefschetz fibrations given by Donaldson's
construction, we can fix this $m$ to be as large as we want at the
start of the argument.

\begin{prop} \label{combine}
If $\alpha^2=\kappa\cdot\alpha=0$ and $j$ is $m$--nondegenerate
for each $m$ under consideration, the total contribution of all
disjoint unions of possibly-multiply-covered tori in classes
proportional to $PD(\alpha)$ to the standard surface count
$\mathcal{DS}_{(X,f)}(n\alpha)$ is the coefficient of $z^n$ in the
product \[ \prod_k
\prod_{C\in\mathcal{M}^{j,\varnothing}_{X}(k\alpha)}P'_j (C,z^k).
\] \end{prop}
\begin{proof} Let $C_i$ be $j$--holomorphic tori in class $k_i
\alpha$, and write $r=\alpha\cdot (fiber)$.  The contribution of a
disjoint union of $m_i$--fold covers of the $C_i$ to\break
$\mathcal{DS}_{(X,f)}(\sum m_i k_i \alpha)$ may be found by using
an almost complex structure $J$ on $X_{\sum m_i k_i r}(f)$
obtained by pushing forward generic smooth almost complex
structures $J_i$ on the $X_{m_i k_i r}(f)$ via the ``divisor
addition'' map $\prod S^{m_i k_i r}\Sigma_t\to S^{\sum m_i k_i
r}\Sigma_t$. This is because the $J$--holomorphic sections will
just be fiberwise sums of the $J_i$--holomorphic sections, which
are in turn close to the sections $s_{C_i}$, and the $C_i$ are
assumed disjoint, so that $J$ (which is smooth away from the
diagonal) will be smooth near each $J$--holomorphic section,
putting us in the situation of Proposition \ref{obvious}. We may
then conclude that the total contribution of such a disjoint union
of covers is $\prod_{i}r'_{j}(C,m_i)$, since $J$--holomorphic
sections are obtained precisely by adding together
$J_i$--holomorphic sections under the divisor addition map, and
there are
 $\prod_{i}r'_{j}(C,m_i)$ ways to do this.
Organizing these contributions into a generating function then yields the
proposition.
\end{proof}

We now fix an embedding $u$ of a square-zero torus $C$ and
consider paths $j_t$ ($t\in (-\ep,\ep)$) of almost complex
structures making $u$ and $f$ holomorphic.  If $(C,j_t)$ crosses a
wall at $t=0$ we would like to compare the $r'_{j_t}(C,m)$ for
small negative values of $t$ to those for small positive values.
We note again that we are taking $\nu_t=0$ for ease of exposition,
but the following lemma and its proof go through unchanged to the
case when we instead have a family $(j_t,\nu_t)$ of admissible
pairs with $(C,j_t,\nu_t)$ crossing a wall just at $t=0$.

\begin{lemma} \label{D0} Assume that $(C,j_t)$ crosses the wall
$\mathcal{D}_0$ at $t=0$ and that the path $j_t$ is generic among
paths of almost complex structures making both $C$ and $f$
holomorphic. Then there is a path of $j_t$--holomorphic tori $C_t$
such that:
\begin{itemize} \item[\rm(1)] For each $t$ the set of
$j_t$--holomorphic tori homologous to $C$ in a suitably small
tubular neighborhood $U$ of  $C$ is $\{C,C_t\}$. \item[\rm(2)]
$C_0=C$ \item[\rm(3)] For $0<|t|<\ep$, $(C_t, j_{-t})$ and $(C,j_t)$
are connected by a path $(C'_s,j'_s)$ with every $j'_s$ making $f$
holomorphic and every $C'_s$ $m$--nondegenerate.
\end{itemize}  Moreover, there are small regular perturbations
$j'_{t}$ of the path $j_t$ supported near $t=0$ with the property
that there are no $j'_0$--holomorphic curves in any homology class
$k[C]$ contained in $U$ \end{lemma}

\begin{proof}
We mimic the argument on pp. 863--864 of \cite{Taubes}.  Let $D$
be the linearization of $\dbar_{j_0}$ at the embedding $u$ of $C$.
For small $|t|>0$, the equation for a section $v_t$ of $N_C$ to
have the property that $\exp_u v_t$ is $j_t$--holomorphic has the
form \begin{equation} \label{lin} Dv_t + R(t,v_t,\nabla v_t)=0
\end{equation} where the Taylor expansion of $R$ begins at order 2
(in the case considered in \cite{Taubes} there is an additional
term proportional to $t$ times the derivative with respect to $t$
of the projection to $N_C$ of the restriction of $j_t$ to
$T_{0,1}C$, but in the present context this term vanishes since
all the $j_t$ make $C$ holomorphic.)  Generically $D$ will have a
one-dimensional kernel and cokernel, so let $s$ span $\ker D$ and
write $v_t = a s + w$ where $a$ is small and $w$ is
$L^2$--orthogonal to $s$; the implicit function theorem lets us
solve the equation obtained by projecting (\ref{lin}) orthogonal
to $\cok D$ for $w$ in terms of $t$ and $a$, so to determine the
structure of the $j_t$ moduli space it remains to solve for $a$ in
terms of $t$.  Now when we project (\ref{lin}) onto $\cok D$ we
obtain an identification of the moduli space in question with the
zero set of a function whose Taylor series begins \begin{equation}
c_1 t^2 + c_2 ta +c_3 a^2. \label{Taylor} \end{equation} Now since
$a=0$ is a solution for all $t$ (corresponding to the curve $C$,
which is $j_t$--holomorphic for all $t$), we have $c_1=0$. Since
$(C,j_t)$ is nondegenerate except at $t=0$, the solution $a=0$ is
nondegenerate for $t\neq 0$, which forces $c_2 \neq 0$.  Moreover,
as in \cite{Taubes}, $c_3 \neq 0$ because of the transversality of
the path $j_t$ to the wall.  It follows that provided the tubular
neighborhood $U$ and the interval $(-\ep, \ep)$ are taken small
enough, the $j_t$--moduli space is as described in the statement
of the Lemma.

Moreover, since the two zeros $a$ of $c_2 ta+c_3 a^2$ are
oppositely-oriented, for each $t$ the spectral flows of the
linearizations at $C$ and $C_t$ of $\dbar_{j_t}$ will be opposite.
Since the sign of the spectral flow for $C$ changes as $t$ crosses
zero, the spectral flows of $(C,j_{-t})$ and $(C_t , j_t)$
therefore have the same sign (ie, the number of
eigenvalue crossings that occur in the flow is the same modulo 2).
Now consider the path \begin{equation} \label{path0} t\mapsto
\left\{
\begin{array}{ll} (C,j_t) & t\leq 0 \\ (C_t,j_t) & t\geq 0
\end{array} \right. \end{equation} The only $t$ at which this path touches a
wall is $t=0$, and we know that the signs of the spectral flows at
the endpoints are the same.  Although curves whose spectral flows
have the same sign may in general lie in different chambers, when
this happens they are separated by at least two walls, not one, so
it follows that $(C,j_{-t})$ and $(C_t, j_t)$ must lie in the same
chamber when $0<t<\ep$ (and, by an identical argument, when
$-\ep<t<0$ as well).  An appropriate perturbation of (\ref{path0})
to a path remaining in this chamber will then have the property
stated in part 3 of the lemma.

For the final part of the lemma, consider generic paths
$\tilde{j}_s$ of almost complex structures with $\tilde{j}_0 =
j_0$ but with the other $\tilde{j}_s$ no longer constrained to
make $C$ holomorphic.  Then exactly as in \cite{Taubes} the moduli
space of $\tilde{j}_s$--holomorphic curves near $C$ will be, for
small $s$, diffeomorphic to the zero set of a function of $a$
whose Taylor series begins $r_1 s + r_2 a^2$ where $r_1$ and $r_2$
are nonzero numbers.  Taking the sign of $s$ appropriately, we
obtain arbitrarily small regular perturbations $\tilde{j}$ of
$j_0$ making no curve near $C$ and homologous to $C$ holomorphic.
By taking $U$ small, we can ensure that there were no embedded
$j_0$--holomorphic curves in any class $k[C]$ where $k>1$ meeting
$U$ (this uses the fact that generically $(C,j_0)$ will not be
located on any of the walls $\mathcal{D}_{\iota}$ with $\iota\neq
0$); if the perturbation $\tilde{j}$ of $j_0$ is taken small
enough there will also not be any $\tilde{j}$--holomorphic curves
meeting $U$ in any of these classes.  Taking a generic
perturbation of the path $j_t$ supported close to zero which
passes though $\tilde{j}$ at $t=0$ then gives the desired result.
\end{proof}

\begin{cor} \label{recip} In the context of Lemma \ref{D0},  for $0<|t|<\ep$, \[
P'_{j_{-t}}(C,z)=\frac{1}{P'_{j_{t}}(C,z)} \] \end{cor}

\begin{proof}
By the third statement in Lemma \ref{D0} and by Corollary
\ref{nowalls}, we have $P'_{j_{-t}}(C,z)=P'_{j_t}(C_t,z)$.  Use a
perturbation $j'_s$ on $U$ of the path $j_s$ as in Lemma \ref{D0}
which differs from $j_s$ only for $|s|<t/2$.  Assuming the
perturbation to be small enough, we may extend $j_s$ and $j'_s$
from the tubular neighborhood $U$ to all of $X$ in such a way that
both are regular outside the neighborhood $U$ (for all $s$) and
they agree with each other outside a slightly smaller region $V$
such that no $j_s$-- or $j'_s$--holomorphic curves are contained
in $U\setminus V$. The contributions of all the $j'_s$ holomorphic
curves outside $U$ will then be constant in $s$.  Since we can use
either $j_{-t} =j'_{-t}$ or $j'_0$ to evaluate the invariant
$\mathcal{DS}$, it follows that the contributions of curves
\emph{inside} $U$ will be the same for $j_{-t}$ as for $j'_0$.
Since the former is obtained from the generating function
$P'_{j_{-t}}(C,z)P'_{j_{-t}}(C_t,z)=P'_{j_{-t}}(C,z)P'_{j_{t}}(C,z)$
while the latter is given by the generating function 1 (for there
are no $j_0$ curves in any class $k[C]$ in the region $U$), the
corollary follows. \end{proof}

Let us now recall some more details in the definition of $Gr$ from
\cite{Taubes}.  The multiple covers of a $j$--holomorphic
square-zero torus $C$ are given weights $r_j (C,m)$ which are
determined by the signs of the spectral flows of each of the four
operators $D_{ \iota}$ to a complex linear operator.  Note that
although Taubes did not define a contribution $r_{j,\nu}(C,m)$
when $\nu\neq 0$, these can be defined using the formulas of
\cite{IP}, in which Eleny Ionel and Thomas Parker interpret the
Gromov invariant as a combination of the invariants of \cite{RT}
(which count solutions to the inhomogeneous Cauchy--Riemann
equations). As with $r'$, we organize the $r_{j,\nu}(C,m)$ into a
generating function $P_{j,\nu}(C,z)=\sum_{m\geq 0} r_{j,\nu} (C,m)
z^m$. Assume as we may thanks to Corollary \ref{good AC 2} that
there exists an integrable complex structure $j_0$ on a
neighborhood of $C$ that makes both $f$ and $C$ holomorphic, and
let $(j_t,\nu_t)$ be a path of admissible pairs with $C$ $j_t$--
holomorphic that connects $j_0$ to the nondegenerate almost
complex structure $j=j_1$, such that $(C,j_t,\nu_t)$ is transverse
to all walls and meets at most one wall $\mathcal{D}_{\iota}$ at
any given $t$.  Assume the walls are met at $0<t_1 <\cdots <t_n
<1$. From Taubes' definition of $Gr$ and from Lemma \ref{intcase}
and Corollary \ref{nowalls}, we have
\[ P'_{j_t}(C,z)=P_{j_t}(C,z)=\frac{1}{1-z} \mbox{ for } t< t_1 \] (in the inhomogeneous case this uses the
formulas of \cite{IP}; see the proof of Corollary \ref{agree} for
more on this). We also know that if $(C,j_t)$ crosses
$\mathcal{D}_0$ at $t_0$, then $P$ and $P'$ both satisfy the
transformation rule
\[ P_{j_{t_0 +\ep}}(C,z)=\frac{1}{P_{j_{t_0 -\ep}}(C,z)} \quad
P'_{j_{t_0 +\ep}}(C,z)=\frac{1}{P'_{j_{t_0 -\ep}}(C,z)}. \] So
since $P$ and $P'$ are both unchanged when $(C,j)$ varies within a
chamber, to show that they agree we need only show that they
transform in the same way when $(C,j_t)$ crosses one of the walls
$\mathcal{D}_{\iota}$ where $\iota\neq 0$.  To again make contact
with the inhomogeneous situation, note that just as the
independence of $\mathcal{DS}$ from the almost complex structure
and the perturbation on $X_r(f)$ used to define it lead to the
wall crossing formulas for the $P'_{j,\nu}$, if we view $Gr$ as a
combination of Ruan--Tian invariants, the independence of these
invariants from the almost complex structure and the perturbation
on $X$ can be considered to lead to wall crossing formulas for the
$P_{j,\nu}$ which are identical to the wall crossing formulas
written down by Taubes in the case $\nu=0$.

We now record the following results, which summarize relevant
parts of Lemmas 5.10 and 5.11 of \cite{Taubes} and their proofs.

\begin{lemma} \label{doublecover} Assume that $(C,j_t)$ crosses the wall
$\mathcal{D}_{\iota}$ where $\iota\neq 0$ at $t=t_0$.  For $\ep$
sufficiently small, $|t-t_0|<\ep$, and for a suitably small
neighborhood $U$ of $C$: \begin{itemize} \item[\rm(1)] The only
connected embedded $j_t$--holomorphic curve homologous to $C$ and
meeting $U$ is $C$ itself. \item[\rm(2)] The only connected, embedded
$j_t$--holomorphic curves meeting $U$ in any homology class $m[C]$
where $m>1$ come in a family $\tilde{C}_t$ in class $2[C]$ defined
either only for $t>t_0$ or only for $t<t_0$.  As $t\to t_0$,
suitably chosen embeddings $\tilde{u}_t\co \tilde{C}_t\to X$
converge to $u\circ \pi\co  \tilde{C}_0\to X$, where $u$ is the
embedding of $C$ and $\pi\co  \tilde{C}_0\to C$ is a double cover
classified by $\iota \in H^1 (C,\mathbb{Z}/2)$. \item[\rm(3)] The
signs of the spectral flows for $(\tilde{C}_{t_0+\delta},
j_{t_0+\delta})$ are the same as those for
$(\tilde{C}_{0},j_{t_0-\delta})$, where $\tilde{C}_{0}$ is mapped
to $X$ by $u_{t_0-\delta}\circ\pi$ (here $\delta$ is any small
number having whatever sign is needed for $\tilde{C}_{t_0+\delta}$
to exist). \end{itemize} \end{lemma}

Using the information from part 3 of the above lemma, the $r_j
(C,m)$ are defined in such a way as to ensure that
\begin{equation} \label{icross}
P_{j_{t_0-\delta}}(C,z)=P_{j_{t_0+\delta}}(C,z)P_{j_{t_0+\delta}}(\tilde{C}_{t_0+\delta},z),\end{equation}
which is necessary for $Gr$ to be independent of the almost
complex structure used to define it; Taubes finds necessary and
sufficient conditions in which the $r_{j}(C,m)$ should depend on
the signs of the spectral flows in order for (\ref{icross}) to
hold.  Meanwhile, the fact that $\mathcal{DS}$ is known a priori
to be independent of the almost complex structure $J$ used to
define it ensures that \begin{equation} \label{i'cross}
P'_{j_{t_0-\delta}}(C,z)=P'_{j_{t_0+\delta}}(C,z)P'_{j_{t_0+\delta}}(\tilde{C}_{t_0+\delta},z),\end{equation}
as can be seen by the usual method of taking smooth almost complex
structures $J_t$ which are H\"older-close enough to the
$\mathbb{J}_{j_t}$ that a $J_t$--holomorphic section in the
relevant homotopy classes meets the neighborhood $\mathbb{U}$ if
and only if it contributes to one of the terms in (\ref{i'cross}),
in which case it is contained in $\mathbb{U}$.  If we somehow knew
a priori that the $r'_j (C,m)$ depended only on the signs of the
spectral flows, then because Taubes' conditions are
\emph{necessary} in order to get an invariant it would follow that
$P'_{j_t}(C,z)$ has to change as $t$ crosses $t_0$ in the same way
that $P_{j_t}(C,z)$ changes.  However, we only know that the
$r'_{j}(C,m)$ are unchanged if we move $(C,j)$ within a chamber;
nonetheless it's not difficult to push what we know far enough to
get the right transformation rule.

\begin{lemma} \label{zsquared} In the context of Lemma \ref{doublecover}, \[
P'_{j_{t_0+\delta}}(C,z)=\frac{P'_{j_{t_0-\delta}}(C,z)}{P'_{j_{t_0-\delta}}(C,z^2)}.\]
\end{lemma}

\begin{proof}
Assume that $(C,j_t)$ crosses some $\mathcal{D}_{\iota}$ with
$\iota\neq 0$ precisely at the point $t_0$, and work in the
notation of Lemma \ref{doublecover}. Observe that, analogously to
the situation for crossings of $\mathcal{D}_0$, since (where
$\delta$ is small and of whichever sign is necessary for the
following statements to make sense)
$(\tilde{C}_{t_{0}+\delta},j_{t_{0}+\delta})$ and
$(\tilde{C}_0,j_{t_{0}-\delta})$ have identical signs for their
spectral flows, and since the path \begin{equation} \label{path}
t\mapsto \left\{
\begin{array}{ll} (\tilde{C}_0,j_t) & t\mbox{ between  } t_{0}-\delta \mbox { and } t_{0} \\
(\tilde{C}_t,j_t) & t\mbox{ between  } t_{0} \mbox { and }
t_{0}+\delta
\end{array} \right. \end{equation}
only meets a wall at $t=t_0$,
$(\tilde{C}_{t_{0}+\delta},j_{t_{0}+\delta})$ and
$(\tilde{C}_0,j_{t_{0}-\delta})$ must lie in the same chamber
(their having identical signs for their spectral flows but lying
in different chambers would require any path between them to meet
two walls).  We can therefore perturb the path (\ref{path}) near
$t_0$ to one (say $t\mapsto (C'_t,j'_t)$)  which stays entirely
within that chamber, with each $j'_t$ making the restriction of
$f$ to the neighborhood of $C'_t$ on which it is defined
pseudoholomorphic.  Hence by Corollary \ref{nowalls} we have
$r'_{j_{t_0+\delta}}
(\tilde{C}_{t_0+\delta},m)=r'_{j_{t_0-\delta}}(\tilde{C}_0,m)$.
But $\tilde{C}_0$ is a double cover of $C$, so in fact
$r'_{j_{t_0+\delta}}
(\tilde{C}_{t_0+\delta},m)=r'_{j_{t_0-\delta}}(C,2m)$,
ie,
\[ P'_{j_{t_0+\delta}}(\tilde{C}_{t_0+\delta},z)=P'_{j_{t_0-\delta}}(C,z^2).\] The
lemma then follows immediately from equation \ref{i'cross}.
\end{proof}

Again, the same wall crossing formula for the $P_{j,\nu}$ for
general $\nu$ follows in exactly the same way, using the
independence of $Gr$ from the data used to define it via the
``Ruan--Tian series'' that appears in \cite{IP}.

\begin{cor} \label{agree} Let $j$ be an almost complex structure as in
Corollary \ref{good AC 2} and $C$ a $j$--holomorphic square-zero
torus.  Then $r'_{j,\nu} (C,m)=r_{j,\nu} (C,m)$ for all $m$ and
$\nu$ for which $(j,\nu)$ is admissible and $(C,j,\nu)$ is
$m$--nondegenerate.
\end{cor}

\begin{proof}  Let $j_t$ be a path of almost complex structures
making $f$ and $C$ holomorphic beginning at an almost complex
structure $j_0$ which is integrable near $C$ and ending at
$j=j_1$, and let $\nu_t$ be inhomogeneous terms such that each
$(j_t,\nu_t)$ is admissible and $(C,j_t,\nu_t)$ is transverse to
all walls; Lemmas \ref{intexists} and \ref{paths} ensure the
existence of such paths.  Assume the walls are crossed at the
points $t_1<\cdots<t_n$ (so that in particular $(C,j_0,\nu_0)$ is
$m$--nondegenerate). Now it follows from the description of the
Gromov invariant in terms of the Ruan--Tian invariants in
\cite{IP} that $r_{j_0,\nu_0}(C,m)=1$ for all $m$: Definition 3.3
and Theorem 4.5 of that paper show that the contribution in
question may be computed by assigning to the various $m$--fold
covers of $C$ (including the disconnected ones) weights which add
up to 1 when all the linearizations of the inhomogeneous equations
are surjective and complex linear. So by Lemma \ref{intcase}, for
all $m$ $P_{j_t,\nu_t}(C,z)=P'_{j_t,\nu_t}(C,z)=\frac{1}{1-z}$ for
all $m$ and all suitably small $t$, and by Corollary \ref{nowalls}
$P_{j_t,\nu_t}(C,z)$ and $P'_{j_t,\nu_t}(C,z)$ change only when
$t$ passes one of the $t_i$. By Corollary \ref{recip} and the
construction of $Gr$ (specifically Equation 5.26 of
\cite{Taubes}), if the wall $\mathcal{D}_0$ is crossed at $t_i$
the changes in both $P$ and $P'$ are found by taking the
reciprocal, while Lemma \ref{zsquared} above and Equation 5.28 of
\cite{Taubes} tell us that if the wall $\mathcal{D}_{\iota}$ with
$\iota\neq 0$ is crossed at $t_i$ then both $P$ and $P'$ change
according to the rule
\[ P_{j_{t_i+\delta},\nu_{t_i+\delta}}(C,z)=\frac{P_{j_{t_i-\delta},\nu_{t_i-\delta}}(C,z)}{P_{j_{t_i-\delta},\nu_{t_i-\delta}}(C,z^2)}
,\] $\delta$ being small and of the same sign as in Lemma
\ref{doublecover}.  Hence
$P'_{j_1,\nu_1}(C,z)=P_{j_1,\nu_1}(C,z)$, proving the corollary.
\end{proof}

The objects which contribute to $Gr(\alpha)$ are, for generic
almost complex structures $j$, formal sums of form $h=\sum m_i
C_i$ where the $C_i$ are disjoint $m_i$--nondegenerate
$j$--holomorphic curves, the $m_i$ are positive integers which are
required to equal 1 unless $C_i$ is a square zero torus, and $\sum
m_i [C_i ]=PD(\alpha)$. For curves $C_i$ which are not square zero
tori, let $r_j(C,1)$ be the contribution of $C$ to $Gr$
(ie, the sign of the spectral flow of the linearization
of $\dbar_j$), and (assuming $j$ makes $f$ holomorphic and
$\mathbb{J}_j$ is regular for $C$) $r'_j (C,1)$ the contribution
of $C$ to $\mathcal{DS}$, so that, by Theorem \ref{singleagree},
$r'_j (C,1)=r_j (C,1)$.  By definition, the contribution of the
formal sum $h$ to $Gr(\alpha)$ is $\prod_i r_j (C_i,m_i)$, while
the proof of Proposition \ref{combine} shows that the contribution
of $h$ to $\mathcal{DS}_{(X,f)}(\alpha)$ is $\prod_i r'_j
(C_i,m_i)$. Thus the previous proposition shows that \emph{every
object $h$ which contributes to $Gr$ contributes to $\mathcal{DS}$
in the same way}. To prove that $\mathcal{DS}=Gr$, we need to see
that, if we compute $\mathcal{DS}$ using an almost complex
structure $J$ H\"older close to a generic $\mathbb{J}_j$, then the
only sections contributing to $\mathcal{DS}$ may be viewed as
contributions from some disjoint union of $j$--holomorphic curves
in $X$ with only square-zero tori allowed to be multiply covered.

To see this, note that for any $\alpha\in H^2 (X,\mathbb{Z})$, by
Gromov compactness, if $J$ is close enough to $\mathbb{J}_j$ then
any $J$--holomorphic sections in the class $c_{\alpha}$ must be
contained in some small neighborhood of a section which
tautologically corresponds to some (generally disconnected, not
embedded) curve in $X$ with total homology class $PD(\alpha)$. Now
for generic $j$, the space of (possibly disconnected)
$j$--holomorphic curves in $X$ which have any singularities
(including intersection points of different connected components)
or have any components other than square-zero tori or exceptional
spheres which are multiply covered has dimension strictly less
than the dimension $d(\alpha)$ (This follows by easy algebra using
the formula for $d(\alpha)$, and is of course the reason that $Gr$
is not obliged to count singular curves or multiply-covered curves
other than square-zero tori). Curves in $X$ with multiply-covered
exceptional sphere components may similarly  be eliminated by a
dimension count: If $\alpha$ is any class represented by a
$j$--holomorphic curve and $\beta$ is the class of an exceptional
sphere, we have $d(\alpha-m\beta)=d(\alpha)-\beta\cdot m\alpha
-\frac{1}{2}(m^2+m)<d(\alpha)-1$, so for generic choices of
$d(\alpha)$ points in $X$, no union $C$ of a $j$--holomorphic
curve in class $\alpha-m\beta$ with an $m$--fold cover of the
$j$--holomorphic sphere in class $\beta$ passes through all
$d(\alpha)$ of the points.

 Hence in any case, the
space of $\mathbb{J}_j$--holomorphic sections tautologically
corresponding to curves not counted by $Gr$ has dimension less
than the dimension of the space of sections counted by
$\mathcal{DS}_{(X,f)}(\alpha)$, which is equal to $d(\alpha)$ by
Proposition 4.3 of \cite{Smith}.  In principle, it perhaps could
happen that when we perturb $\mathbb{J}_j$ to a smooth almost
complex structure $J$ near such a section $s_C$ to find the
contribution of $C$ we might obtain a positive-dimensional set of
nearby $J$--holomorphic sections, but because these sections are
constrained by Gromov compactness to stay near $s_C$, for a large
open set of choices of the incidence conditions used to cut down
the moduli spaces for $Gr$ and $\mathcal{DS}$ to be
zero-dimensional, the perturbed sections will still not appear in
this moduli space and so will not contribute to $\mathcal{DS}$.

$\mathcal{DS}$ and $Gr$ therefore receive contributions from just
the same objects, so since these contributions are equal, Theorem
\ref{main} follows.

\end{document}